\documentclass{article}
\date{}

\usepackage{graphicx}

\newcommand{\bey}{\begin{eqnarray}}
\newcommand{\eey}{\end{eqnarray}}
\newcommand{\nn}{\nonumber}
\newcommand{\beq}{\begin{equation}}
\newcommand{\eeq}{\end{equation}}
\newtheorem{thm}{\hspace{6mm}Theorem}[section]

\newtheorem{lem}{\hspace{6mm}Lemma}[section]

\newcommand{\proofend}{\mbox{ }\hfill \raisebox{.4ex}{\framebox[1ex]{}}}

\def\e{{\epsilon}}

\def\no{{\nonumber}}
\def\a{{\alpha}}

\def\s{{\sigma}}

\def\min{\mbox{min}}

\title{A posteriori error analysis for finite element solution of
elliptic differential equations using equidistributing meshes}

\author{Yinnian He
\thanks{Faculty of Science, Xi'an Jiaotong
University, Xi'an 710049, People's Republic of China
(heyn@mail.xjtu.edu.cn).}
\and Weizhang Huang
\thanks{Department of Mathematics, University of Kansas, Lawrence, KS 66049, U.S.A. (huang@math.ku.edu).}}

\begin{document}
\vskip 1cm
\maketitle

\begin{abstract}
The paper is concerned with the adaptive finite element solution of linear elliptic differential equations
using equidistributing meshes. A strategy is developed for defining this type of mesh
based on residual-based a posteriori error estimates and rigorously analyzing
the convergence of a linear finite element approximation using them.
The existence and computation of equidistributing meshes and
the continuous dependence of the finite element approximation on mesh are also studied.
Numerical results are given to verify the theoretical findings.
\end{abstract}

\noindent{\bf AMS 2000 Mathematics Subject Classification.}
65M50,65M60, 65L50, 65L60

\noindent{\bf Key Words.}
Mesh adaptation, equidistribution, error analysis, finite element method

\noindent{\bf Abbreviated title.}
Error analysis for equidistributing meshes

%----------------------- section 1 ---------------------------
\section{Introduction}
\label{SEC:introduction}

We are concerned with the convergence of the linear finite element
solution of elliptic differential equations using equidistributing
meshes. An equidistributing mesh
of $N$ elements for $\Omega \equiv (0,1)$ is a mesh $x_0 = 0 < x_1 < \cdots < x_N = 1$
satisfying the so-called equidistribution principle \cite{Bur74,deB73}
\beq
\int_{x_{i-1}}^{x_i} \rho(x) d x =
\frac{1}{N}\int_0^1 \rho(x) d x, \quad i = 1, ..., N
\label{eq-1}
\eeq
where $\rho =\rho(x)$ is a user-prescribed, strictly positive
function. Function $\rho(x)$, referred to as an adaptation function,
can be interpreted as an ``error'' density function, with
$\int_0^1 \rho(x) d x$ being the total ``error''. Equation
(\ref{eq-1}) implies that $\rho(x)$ is evenly distributed among the mesh elements.

Equidistributing meshes are known to produce optimal error bounds and have been
widely used for adaptive numerical solution of differential equations.
Their theoretical studies have also attracted considerable attention from researchers;
e.g., see \cite{Bur74,deB73,deB74,Dod72,Mcc70,Ric69} on best approximations
with variable nodes, \cite{SY66,SY68,SY70} on regression problems in statistics,
\cite{BR79b,PS75,Whi79} on adaptive numerical solution of differential equations,
and  \cite{BM00,BM01,BM01b,CSX06,Hua03,Hua02b,HS01,KS01,Mac99,QS99,QST00}
for more recent works.
A focus of these studies has been on error analysis, i.e., 
to understand how accurate an approximation or a numerical solution can be
on an equidistributing mesh.
Unfortunately, this has proven to be a difficult task due to the highly
nonlinear coupling between the mesh and the solution. 
The analysis can be significantly simplified by taking
a priori meshes defined using the exact solution or some information of the exact solution.
Interestingly, almost all of the existing analyses have been done in this way.
For example, Pereyra and Sewell \cite{PS75} choose a mesh to equidistribute 
a form of the truncation error and obtain an asymptotical bound for it
for the finite difference solution of two-point boundary value problems.
Qiu et al. \cite{QS99,QST00} and Beckett and Mackenzie \cite{BM00,BM01,BM01b,Mac99}
investigate the uniform convergence of finite difference and finite element approximations for
singularly perturbed problems for meshes determined using the equidistribution principle
and the singular part of the exact solution.
Chen and Xu \cite{CX08} show that a standard finite element method and a new streamline
diffusion finite element method produce stable and accurate approximations
for a singularly perturbed convection-diffusion problem provided that
the mesh properly adapts to the singularity of the solution.
Huang et al. \cite{Hua03,Hua02b,HS01} and Chen et al. \cite{CSX06}
study multi-dimensional interpolation problems using equidistributing meshes which depend on
the function under consideration.

The noticeable exceptions are the work \cite{BR79b} and \cite{KS01} where
a posteriori equidistributing meshes, or equidistributing meshes determined by the computed solution,
are considered. More specifically, Babu\u{s}ka and Rheinboldt \cite{BR79b}
consider the linear finite element solution
of a one-dimensional elliptic problem and develop a functional
from a residual-based a posteriori error estimate in lieu of asymptotic approximation and
coordinate transformation. Using the optimal coordinate transformation obtained by minimizing
the functional,
Babu\u{s}ka and Rheinboldt show that a mesh is asymptotically optimal
if the residual-based error estimate is evenly distributed among the mesh elements.
Kopteva and Stynes \cite{KS01} study an upwind finite difference discretization
of  one-dimensional quasi-linear convection-diffusion problems without turning points
and develop a convergence analysis for the discretization
where the mesh is determined by the computed solution through the equidistribution
principle and the arc-length adaptation function.

In this paper we are concerned with convergence analysis for the finite element solution using
a posteriori equidistributing meshes. The goal is to develop a systematic approach for
defining these meshes such that both their error analysis and
computation can be done in an a posteriori manner.
At the same time we would like the approach to be general enough so that
it can apply to any standard finite element method and have no essential limitations
for multi-dimensional generalizations. Furthermore, the approach should be
mathematically rigorous. Particularly, it should not rely on
asymptotic approximation or continuous coordinate transformations as in \cite{BR79b,HS01}.
Several  other issues, such as the existence and computation of equidistributing meshes
and the continuous dependence of the linear finite element solution on mesh, are also
studied in the paper. The main results are given in \S\ref{SEC:main-results}.

Since D\"orfler's seminal work \cite{Dor96} significant progress has been made 
on the convergence analysis of adaptive finite element methods based on a posteriori error estimates;
e.g. see \cite{BDD04,CH06,CKNS08,CHX08,MNS00,MNS02,Ste05}. It should, however, be pointed out that
there are essential differences between those works and the current one. 
The former ones are dealt with adaptive mesh refinement using specially designed
marking strategies and their convergence results are typically measured in terms of refinement levels,
whereas the current work is concerned with equidistributing meshes
(including their existence, generation, optimality, and error analysis) and our results are measured in terms of
the number of mesh elements (cf. Theorems \ref{thm2.1} and \ref{thm2.2}).
It does not seem that the existing convergence analysis for mesh refinement
can apply directly to equidistributing meshes and neither can the current results be covered by
the existing ones. On the other hand, adaptive mesh refinement and equidistribution do share some
common ground. For example, an equidistributing mesh can be generated through
mesh refinement (e.g., see \cite{BGHLS97,Hua02b}) (and other strategies (e.g. see \cite{Hua01b}
for a variational approach)), and the concept of mesh equidistribution is often used in mesh refinement
algorithms and computer codes for maximizing the efficiency of computation (e.g., see \cite{MA91}).
Relations between convergence results for adaptive mesh refinement and equidistribution may thus
deserve further investigations.

The paper is organized as follows. The description of the mathematical problem and the main results
are given in \S\ref{SEC:main-results}. The approach for defining equidistributing meshes and
analyzing the corresponding finite element error is developed in \S\ref{SEC:error-analysis}.
An iterative algorithm for computing the meshes is proposed and numerical results are
presented in \S\ref{SEC:numerical-results}. The continuous dependence of the finite element
solution on mesh and the existence of equidistributing meshes are studied
in \S\ref{SEC:continuity} and \S\ref{SEC:existence}, respectively.

% section 2
\section{Main results}
\label{SEC:main-results}

We consider the boundary value problem of a linear elliptic
differential equation
\bey && - (a u')' + b u' + c u = f,\qquad \mbox{ in } \Omega
\label{pde-1}
\\
&& u(0) = u(1) = 0 ,
\label{bc-1}
\eey
where $a(x)$, $b(x)$, $c(x)$, and
$f(x)$ are given functions satisfying
\beq \quad a, b \in
W^{1,\infty}(\Omega), \quad c\in L^{\infty}(\Omega),\quad f \in
L^2(\Omega), \label{coef-1} \eeq and \beq a(x) \ge a_0 > 0, \quad
c(x) - \frac{1}{2} b'(x) \ge 0, \quad a.e.\mbox{ in } \Omega
\label{coef-2} \eeq for some constant $a_0$. Here,
$W^{1,\infty}(\Omega)$ denotes the Sobolev space of functions whose
derivatives are in $L^\infty(\Omega)$. The variational form of
problem (\ref{pde-1}) and (\ref{bc-1}) is to find $u\in V \equiv
H^1_0(\Omega)$ such that \beq B(u,v)=(f,v), \quad \forall v\in V
\label{bvp-1} \eeq where \beq B(u,v) = \int_\Omega (a u' v' + b u'v
+ cuv ) d x, \quad (f, v) = \int_\Omega f v d x. \label{bform-1}
\eeq For a given mesh \beq \pi_h:\quad x_0 = 0 < x_1 < \cdots < x_N
= 1 \label{mesh-1} \eeq with $h=\max_i ( x_i-x_{i-1})$,
a linear finite element approximation to the solution of
(\ref{bvp-1}) is defined as $u_h \in V_h$ satisfying \beq
B(u_h,v_h)=(f,v_h), \quad \forall v_h\in V_h \label{fem-1} \eeq
where the linear finite element space $V_h$ is given by \beq V_h =
\mbox{span}\{\phi_1, ..., \phi_{N-1}\} , \label{Vh-1} \eeq with
$\phi_i$'s being linear basis functions associated with mesh points
$x_i$'s.

We are concerned with adaptive finite element solution of
(\ref{pde-1}) and (\ref{bc-1}) using equidistributing meshes. For
this purpose, we choose the mesh according to the equidistribution
principle (\ref{eq-1}) or
\beq
\rho_{i} h_i = \int_{x_{i-1}}^{x_i}
\rho(x) d x = \frac{\sigma_h}{N}, \quad i = 1, ..., N \label{eq-2}
\eeq
where
\bey && h_i = x_{i} - x_{i-1},
\label{mesh-2}
\\
&& \sigma_h = \sum_i \rho_{i} h_i ,
\label{sigma-1}
\\
&&  \left . \rho \right
|_{(x_{i-1},x_{i})} = \rho_{i} \equiv  \left (1+\frac{1}{\alpha_h }
\left <r_h\right >_i^2\right )^{\frac 1 3},\quad i = 1, ..., N
\label{rho-1}
\\
&& \alpha_h = \left [  \sum_i h_i \left <r_h\right >_i^{\frac 2 3} \right]^3,
\label{alpha-1}
\\
&& r_h  = f + a' u_h' - b u_h'-c u_h ,
\label{res-1}
\\
&& \left <r_h\right >_i=\left (\frac1{h_i}\int^{x_i}_{x_{i-1}}|r_h|^2dx\right )^{\frac 1 2} .
\label{average-1}
\eey
Here, $u_h$ is the solution of (\ref{fem-1}), $r_h$ is the residual function,
$\rho(x)$ is the piecewise constant adaptation function, and
$\left <r_h\right >_i$ is the $L^2$ average of $r_h$ over $(x_{i-1},x_{i})$.

\begin{figure}[thb]
\centering
\includegraphics[width=3.5in]{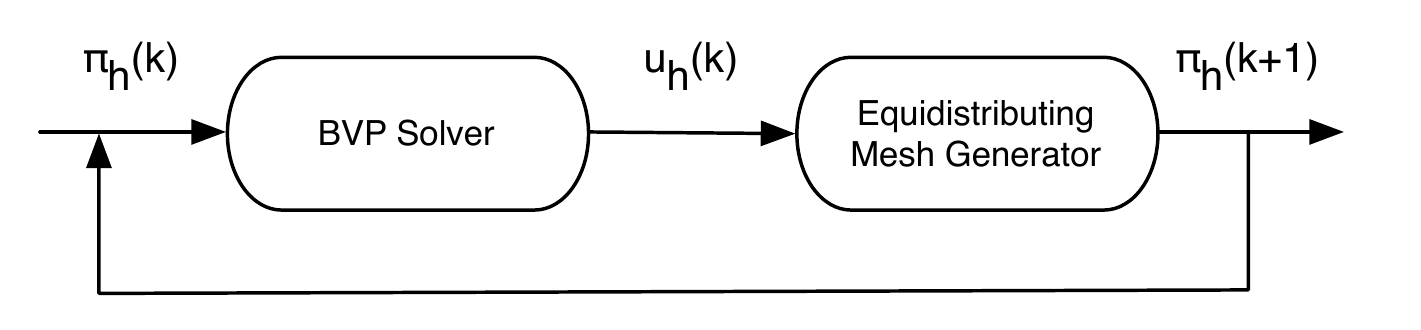}
\caption{Illustration of the iterative solution procedure for the finite
element solution using equidistributing meshes.}
\label{f1}
\end{figure}

The choice (\ref{rho-1}) for the adaptation function is based on an a posteriori
error estimate for the linear finite element solution; see \S\ref{SEC:error-analysis}.
Clearly, the choice depends on the computed solution and thus the mesh is a posteriori.
Unfortunately, this also means that
the mesh and the computed solution are coupled with
each other. The system for $u_h$ and $\pi_h$ consists of algebraic equations
(\ref{fem-1}) and (\ref{eq-2}) and the boundary conditions $x_0 = 0$ and
$x_N = 1$ and is typically solved iteratively; see Fig. \ref{f1}.
An algorithm of this type is given in \S\ref{SEC:numerical-results}.
The existence of the equidistributing mesh is stated in Theorem \ref{thm2.4} below.

% f1

{\em Notation.} The $L^2$-norm on $\Omega$ is denoted by $\| \cdot \|_{\Omega}$
and other $L^p$-norm by $\| \cdot \|_{L^p(\Omega)}$, with the latter being extended to the situation
$0 < p < 1$. Let
\beq
\bar{v } = \| v\|_{L^\infty(\Omega)},\quad \forall v \in L^\infty(\Omega) .
\label{coef-3}
\eeq
We use $C$ as a generic constant which may have different values at different appearances.
In most part of this paper, constants are considered as numbers that may depend
on the domain and coefficients $a(x)$, $b(x)$, and $c(x)$ of differential equation (\ref{pde-1})
but not on the solution $u$, the right-hand sider $f$, and the mesh employed in the finite element
solution. The exceptions are Theorems \ref{thm2.3} and \ref{thm2.4} and
\S\ref{SEC:continuity} and \S\ref{SEC:existence} where constants may further depend on $u$ and $f$.

Define
\beq
\rho_0 =
\left [1+ \gamma \left ( \frac{\|r\|_{L^\infty(\Omega)}}{\|r\|_{L^{\frac23}(\Omega)}}\right )^{\frac23} \right ]^2,
\label{rho0}
\eeq
where $\gamma$ is a positive constant dependent on the domain and coefficients of 
equation (\ref{pde-1}).
The definition of $\gamma$ is given in the proof of Lemma \ref{lem6.0}. The same lemma also shows
that $\rho_0$ is the upper bound on the adaptation function defined in (\ref{rho-1}).

From (\ref{mesh-1}), the mesh $\pi_h$ corresponds one-to-one to the $(N-1)$-component vector
$X = (x_1, ..., x_{N-1})^{\top}$. For this reason,  occasionally $X$ is directly referred to as  a mesh. Let
\bey
S_N = && \left \{ \frac{}{}\right . X \in \Re^{N-1}: \; x_0 = 0 < x_1 < \cdots < x_{N-1} < x_N = 1,
\nn \\
& & \quad \left. \frac{1}{\rho_0 N} \le h_i \le \frac{2}{N},\;\; i = 1, ..., N \; \right \} .
\label{SN}
\eey
It is easy to verify that $S_N$ is a closed, convex subset of $\Re^{N-1}$. The set is equipped with the maximum norm, viz.,
\beq
\| X \|_\infty = \max_i | x_i|, \quad \forall X \in S_N.
\label{SN-norm}
\eeq
This set plays an important role in the study of the existence of equidistributing meshes
and the continuous dependence of the finite element solution on mesh.
Any $N$-cell equidistributing mesh is a member of this set.

\vspace{10pt}

The main results of this paper are summarized in the following four theorems.

\begin{thm}
\label{thm2.1}
{\em (Convergence for equidistributing meshes)}
Define $\rho$ and $\alpha_h$ as in (\ref{rho-1})
and (\ref{alpha-1}), respectively. For any equidistributing
mesh satisfying (\ref{eq-2}), the error for the linear finite element solution (\ref{fem-1})
is bounded by
\beq
\| (u-u_h)'\|_\Omega \le \frac{C \sqrt{\alpha_h} }{N}  ,
\label{thm2.1-1}
\eeq
where $\alpha_h$ has the property
\beq
\lim_{N \to \infty} \sqrt{\alpha_h} = \| r\|_{L^{\frac 2 3}(\Omega)}
\label{thm2.1-2}
\eeq
and  $r$ is the continuous ``residual'' function $r$ defined as
\beq
r = f + a' u' - b u' - c u .
\label{res-2}
\eeq

If further $r$ satisfies $r' \in L^1(\Omega)$, then there exists a positive constant $c$ such that
for $N > c$,
\beq
\left (1+ \left (\frac{c}{N}\right )^{\frac 2 3} \right)^{-\frac{3}{2}}
\| r \|_{L^{\frac 2 3}(\Omega)} \le \sqrt{\alpha_h}
\le \left (1- \left (\frac{c}{N}\right )^{\frac 2 3} \right )^{-\frac{3}{2}} \left [ \| r \|_{L^{\frac 2 3}(\Omega)}^{\frac 2 3}
+ \left (\frac{c \| r' \|_{L^1(\Omega)}}{N}\right )^{\frac 2 3} \right ]^{\frac 3 2} .
\label{thm2.1-3}
\eeq
\end{thm}

\vspace{10pt}

The proof of this theorem is given in \S\ref{SEC:convergence}.
The theorem shows that the error has the asymptotic bound as
\beq
\lim_{N \to \infty} N\; \| (u-u_h)'\|_\Omega \le C \| r \|_{L^{\frac 2 3}(\Omega)}.
\label{thm2.1-10}
\eeq
This is compared with the error bound for a uniform mesh (cf. Lemma \ref{lem3.6})
\beq
\|(u-u_h)'\|_\Omega \le \frac{C_1 \|r \|_\Omega}{ N} .
\label{uni-err-1}
\eeq
Since
\[
\| r \|_{L^{\frac 2 3}(\Omega)} \le \| r \|_{\Omega}
\]
and particularly, the left-hand side is much smaller than the right-hand side when
$r$ ($\sim a u''$) is non-smooth, the theorem implies that the error bound
for an equidistributing mesh can be much smaller than that for a uniform mesh.
This explains why an adaptive mesh often produces a more accurate solution
than a uniform one when the solution is non-smooth.

In practice, it is more realistic to use a quasi-equidistributing
mesh than an exact one. A  quasi-equidistributing mesh is a mesh satisfying
\beq
\frac{N\rho_{i} h_i}{\sigma_h} \le \kappa \qquad i = 1, ..., N
\label{eq-3}
\eeq
for some positive constant $\kappa$
independent of $i$ and $N$. The following theorem, proved in \S\ref{SEC:convergence},
shows that for small $\kappa$,
a quasi-equidistributing mesh leads to a comparable error bound
as an exact equidistributing mesh.

\begin{thm} 
\label{thm2.2} {\em (Convergence for quasi-equidistributing meshes)}
Define $\rho$ and $\alpha_h$ as in (\ref{rho-1})
and (\ref{alpha-1}), respectively. Then for any
quasi-equidistributing mesh satisfying (\ref{eq-3}), the error for the linear finite element solution (\ref{fem-1})
is bounded by
\beq \|
(u-u_h)'\|_\Omega
\le \frac{ C \sqrt{\alpha_h \kappa^3} }{N},
\label{thm2.2-1}
\eeq
where $\alpha_h$ satisfies (\ref{thm2.1-2}).

If further $r$ satisfies $r' \in L^1(\Omega)$, then there exists a positive constant $c$ such that
for $N > \kappa c$,
\beq
\left (1+ \left (\frac{\kappa c }{N}\right )^{\frac 2 3} \right )^{-\frac{3}{2}}
\| r \|_{L^{\frac 2 3}(\Omega)} \le \sqrt{\alpha_h} \le \left (1-
\left (\frac{\kappa c }{N}\right )^{\frac 2 3} \right)^{-\frac{3}{2}} \left [ \| r \|_{L^{\frac 2 3}(\Omega)}^{\frac 2 3}
+ \left (\frac{\kappa c \| r' \|_{L^1(\Omega)} }{N}\right )^{\frac 2 3} 
 \right ]^{\frac 3 2},
\label{thm2.2-3}
\eeq
where $r$ is defined in (\ref{res-2}).
\end{thm}

\begin{thm} {\em (Continuous dependence of finite element solution on mesh)}
\label{thm2.3}
Assume that  $f \in L^\infty(\Omega)$ and $u \in H^2(\Omega)$. Then for any
meshes $X,\; \tilde{X} \in S_N$ satisfying
\beq
\| X - \tilde{X}\|_\infty < \frac{1}{\rho_0 N},
\label{thm2.3-1}
\eeq
the corresponding linear finite element solutions, $u_h$ and $u_{\tilde{h}}$, satisfy
 \bey
\|(u_h-u_{\tilde{h}})'\|_{L^1(\Omega)}&\le&
CN\|X-\tilde{X}\|_\infty,
\label{thm2.3-2}
\\
\|u_h-u_{\tilde{h}}\|_\Omega&\le&
CN^{-\frac12}\|X-\tilde{X}\|^{\frac12}_\infty,
\label{thm2.3-3}
\\
\left | \|u_h\|_E^2-\|u_{\tilde{h}}\|_E^2\right |&\le&
CN^{-\frac12}\|X-\tilde{X}\|^{\frac12}_\infty,
\label{thm2.3-4}
\eey
where $\|\cdot \|_E$ denotes the energy norm associated with the bilinear form $B(\cdot, \cdot)$, viz.,
\[
\| v \|_E^2=\int_\Omega\left ( av'^2+(c-\frac12b')v^2\right )dx .
\]
\end{thm}

\begin{thm} {\em (Existence of equidistributing meshes)}
\label{thm2.4}
Assume that $f\in L^\infty(\Omega)$ and $u \in H^2(\Omega)$. For sufficiently large $N$ (i.e., $N \ge N_0$
where $N_0$ is defined in Lemma \ref{lem6.0}), there exists at least an equidistributing mesh
satisfying (\ref{eq-2}).
\end{thm}

\vspace{10pt}

The above two theorems are proven
in \S\ref{SEC:continuity} and \S\ref{SEC:existence}, respectively.
An iterative algorithm for computing equidistributing meshes is proposed in \S\ref{SEC:numerical-results}.
The numerical results presented in \S\ref{SEC:numerical-results}
demonstrate that the algorithm converges for sufficiently large $N$ and faster for larger $N$.
This is consistent with what observed by Pryce \cite{Pry89} and Xu et al. \cite{XHRW08}
on the convergence of de Boor's algorithm for generating equidistributing meshes
for given adaptation functions.

% section 3
\section{Error analysis for finite element solution using equidistributing meshes}
\label{SEC:error-analysis}

In this section we present an error analysis for equidistributing
meshes satisfying (\ref{eq-2}) and quasi-equidistributing
meshes satisfying (\ref{eq-3}). The approach we use consists of
three major steps, deriving a residual-based a posteriori error
estimate, defining the adaptation function (\ref{rho-1}) based
on the estimate, and developing the error analysis for the
corresponding equidistributing mesh. This approach shares some similarity
with that used in \cite{Hua03,Hua02b,HS01} for analyzing interpolation
error in multi-dimensions. The main difference lies in that the
current analysis is based on an a posteriori error estimate and is
mathematically rigorous, whereas the analysis in
\cite{Hua03,Hua02b,HS01} is based on interpolation error bounds
(depending on the exact solution) and valid only in an asymptotic sense.

% section 3.1
\subsection{Preliminary results}

For completeness and for easy reference we list here some preliminary results without giving their proofs.
These results can be found
in most finite element textbooks, e.g. \cite{BS94,Cia78}.

\begin{lem}
\label{lem3.2} The bilinear form $B(\cdot,\cdot)$ defined in
(\ref{bform-1}) has the properties
\bey && a_0 \, \|v'
\|_{\Omega}^2\leq B(v,v)\le C \| v' \|_\Omega^2, \qquad \forall v
\in V \label{lem3.2-1}
\\
&& B(u,v)\leq C \|u'\|_\Omega\, \| v' \|_\Omega , \qquad
\forall u, v \in V .
\label{lem3.2-2}
\eey
Moreover, the solution of  the problem (\ref{bvp-1}) satisfies
\beq \| u'\|_\Omega \le C \|f\|_{\Omega}.
\label{lem3.2-3}
\eeq
\end{lem}

\begin{lem}
\label{lem3.3}
Given a mesh $\pi_h$, denote by $\Pi_h$ the operator for piecewise linear interpolation, i.e.,
\beq
\Pi_h v (x)  = \sum_{i=0}^N v(x_i) \phi_i(x),
\qquad \forall v \in H^1(\Omega) 
\label{interp-1}
\eeq
where
$\phi_i$'s are the linear basis functions associated with mesh points $x_i$'s.
Then, for any $K_i = (x_{i-1},x_i)$,
\bey
\|v-\Pi_h v\|_{K_i} &
\le &  C h_i \|(v-\Pi_hv)'\|_{K_i}, \qquad  \forall
v\in H^1(K_i) \label{lem3.3-1}
\\
\| (v-\Pi_h v)'\|_{K_i} & \le & C h_i \|v''\|_{K_i},
\qquad \forall v\in H^2(K_i) \label{lem3.3-2}
\\
 \| (v-\Pi_h v)'\|_{K_i} & \le & C \|v'\|_{K_i},
\qquad  \forall v\in H^1(K_i). \label{lem3.3-3} \eey
\end{lem}

\vspace{5pt}

The error for the finite element solution $u_h$, $e_h = u-u_h$, satisfies the orthogonality property and
the error equation, viz.,
\bey
&& B(e_h, v_h) = 0, \qquad \forall v_h \in V_h 
 \label{orth-1}
 \\
&& B(e_h, v) = (f, v) - B(u_h, v), \qquad \forall v \in V.
\label{erreqn-1}
\eey

\vspace{5pt}

\begin{lem}
\label{lem3.6}
The finite element solution $u_h$ defined in (\ref{fem-1}) satisfies
\beq
\| u_h' \|_\Omega \le C \|f\|_{\Omega} .
\label{lem3.6-0}
\eeq
Moreover, if the solution of the continuous problem (\ref{bvp-1}) satisfies $u\in H^2(\Omega)$ and
the mesh has the property
\beq
h \equiv \max_i h_i \le \frac{C_1}{N}
\label{h-1}
\eeq
for some positive constant $C_1$, the error  is bounded by
\beq
\|(u-u_h)'\|_\Omega \le \frac{C \|r\|_\Omega}{N}.
\label{lem3.6-1}
\eeq
\end{lem}

\vspace{5pt}

Obviously, a uniform mesh satisfies the condition (\ref{h-1}).  The finite element
error for a uniform mesh can thus be bounded as in (\ref{lem3.6-1}).

% section 3.2
\subsection{An a posteriori error estimate}

We now derive a residual-based a posteriori error estimate for the
finite element solution. The general procedure for this type of 
error estimation can be seen, e.g.,  in  \cite{AO00,BS01,BW85,BR01,Ver96}.

\begin{lem}
\label{lem3.5} The error $e_h = u-u_h$ is bounded by
 \beq
\| (u-u_h)'\|_\Omega^2 \le C \eta_h^2= C \sum_i h_i^2
\|r_h\|_{K_i}^2, \label{lem3.5-1} \eeq where $r_h$ is defined in
(\ref{res-1}), i.e., $r_h = f + a' u_h' - b u_h' - c u_h$.
\end{lem}

\vspace{5pt}

{\bf Proof.} Using orthogonality property (\ref{fem-1}), error equation (\ref{orth-1}), integration by parts,
Lemma \ref{lem3.3}, and Schwarz' inequality, we have,  for any $v \in V$,
\bey
B(e_h,v)&=&B(e_h,v-\Pi_h v)
\nn \\
& = & (f,v-\Pi_h v)-B(u_h,v-\Pi_h v) 
\nn \\
&=&\sum_{i} \int_{K_i} r_h (v-\Pi_h v)dx
\nn \\
&\le&\sum_{i}\|r_h\|_{K_i} \|v-\Pi_h v\|_{K_i}
\nn \\
&\le&C \sum_{i}\|r_h\|_{K_i} \cdot h_i \|v'\|_{K_i}
\nn \\
&\le& C \left ( \sum_{i}h_i^2 \|r_h\|_{K_i}^2\right
)^{1/2}\| v'\|_\Omega.
\nn
\eey
Then (\ref{lem3.5-1}) follows by taking $v = e_h$ in the above
inequality and using Lemma \ref{lem3.2}.
\proofend

% section 3.3
\subsection{Determination of optimal adaptation function}

Up to this point the mesh has been assumed to be arbitrary and
Lemma \ref{lem3.5} has been obtained for this
general mesh. From now on we shall focus on equidistributing meshes
determined according to the a posteriori error estimate
(\ref{lem3.5-1}).

As we can see from (\ref{eq-1}), the key for the determination of
equidistributing meshes is to define an appropriate adaptation
function $\rho=\rho(x) > 0$. To this end, we regularize $\eta_h$
in (\ref{lem3.5-1}) with a positive constant $\alpha_h > 0$ (to be
determined), i.e.,
\bey \eta_h^2 & = &  \sum_{i}
h_i^2 \|r_h\|_{K_i}^2
\nn \\
& \le & \sum_{i} h_i^3 \left  (\alpha_h +
\left <r_h\right >_i^2\right )
\nn \\
& = & \alpha_h  \sum_{i} h_i^3\left  (1 +
\frac{1}{\alpha_h } \left <r_h\right >_i^2\right ) \equiv \tilde{\eta}_h^2 .
\label{eta-1}
\eey
Notice that $\tilde{\eta}_h \to \eta_h$ as
$\alpha_h \to 0$. Moreover,
\bey
\tilde{\eta}_h^2  & = &
\alpha_h \sum_{i} \left [h_i \left  (1 +
\frac{1}{\alpha_h } \left <r_h\right >_i^2\right )^{\frac 1 3}\right ]^3
\nn \\
& \ge &  \frac{\alpha_h}{N^2} \left [ \sum_{i} h_i \left  (1
+ \frac{1}{\alpha_h} \left <r_h\right >_i^2\right )^{\frac 1 3}\right ]^3 .
\label{eta-2}
\eey
It is easy to see that the equality in the above inequality holds
for any equidistributing mesh when
the adaptation function $\rho(x)$ is chosen as in (\ref{rho-1}), i.e., 
\beq
\left . \rho \right |_{K_i} = \rho_{i}
\equiv \left  (1 + \frac{1}{\alpha_h} \left <r_h\right >_i^2\right
)^{\frac 1 3}, \quad i = 1, ..., N.
\label{rho-2}
\eeq
The bound in (\ref{eta-2}) is not the lowest due to its mesh dependence.
Nevertheless, we may expect
\[
\alpha_h \left [ \sum_{i} h_i \left  (1 + \frac{1}{\alpha_h}
\left <r_h\right >_i^2\right )^{\frac 1 3}\right ]^3 \to \alpha \left [
\int_\Omega\left (1+\frac{1}{\alpha} r^2 \right )^{\frac 1 3} d
x\right ]^3 \quad \mbox{ as } \quad N \to \infty ,
\]
where $\alpha = \lim_{N\to \infty} \alpha_h$. When this is the case,
(\ref{eta-2}) is an asymptotically lowest bound for $\tilde{\eta}_h^2$, and in this sense
the choice (\ref{rho-2}) is asymptotically optimal.

From (\ref{eq-2}), (\ref{eta-2}), and Lemma \ref{lem3.5}, it is easy
to see that the error on a mesh equidistributing the so-defined
$\rho$ is bounded by
\beq
\|(u-u_h)'\|_\Omega^2 \le C \tilde{\eta}_h^2 = \frac{C \alpha_h \sigma_h^3}{ N^2} .
\label{eta-3}
\eeq

To complete the definition, we need to determine the parameter
$\alpha_h$. We follow \cite{Hua01} to choose it such that
\beq
\sigma_h \equiv \sum_i h_i \rho_{i} \le 2.
\label{sigma-2}
\eeq
In this way, roughly fifty percents of the mesh points are placed in the
region where $\rho \gg 1$ \cite{Hua01}. From Jensen's inequality and (\ref{rho-2}),
\bey
\sigma_h & = & \sum_i h_i \left  (1 + \frac{1}{\alpha_h}
\left <r_h\right >_i^2\right )^{\frac 1 3}
\nn  \\
& \le & \sum_i h_i  \left ( 1 + \alpha_h^{-\frac 1 3}
\left <r_h\right >_i^{\frac23}\right )
\nn \\
& = & 1 + \alpha_h^{-\frac 1 3} \sum_i h_i \left <r_h\right >_i^{\frac 2 3}.
\nn
\eey
Then (\ref{sigma-2}) holds when $\alpha_h$ is chosen as in
(\ref{alpha-1}), i.e.,
\beq \alpha_h = \left [ \sum_i h_i
\left <r_h\right >_i^{\frac 2 3}\right ]^3 .
\label{alpha-2}
\eeq

Thus, when the adaptation function $\rho(x)$ and the intensity
parameter $\alpha_h$ are chosen as in (\ref{rho-2}) and
(\ref{alpha-2}), respectively,  from (\ref{eta-3}) and
(\ref{sigma-2}) we see that the finite element error for a mesh equidistributing
$\rho$ is bounded by
\beq \|(u-u_h)'\|_\Omega^2 \le
C \tilde{\eta}_h^2 \le \frac{C \alpha_h}{N^2} .
\label{eta-4}
\eeq
The boundedness of $\alpha_h$ as $N \to \infty$ is investigated in the next subsection.

% section 3.4
\subsection{Convergence for equidistributing and quasi-equidistributing meshes}
\label{SEC:convergence}

We notice that the adaptation function defined in (\ref{rho-2}) 
satisfies $\rho_{i} \ge 1,\, i = 1, ..., N$. As a consequence, (\ref{sigma-2}) implies that
the equidistributing mesh has the property (\ref{h-1}) with $C_1 = 2$.
Combining this with Lemma \ref{lem3.6} we have the following theorem.

\begin{thm}
\label{thm3.1}
Assume that $\rho$ and $\alpha_h$ are defined as in
(\ref{rho-2}) and (\ref{alpha-2}), respectively. If $u \in
H^2(\Omega)$, then for any mesh equidistributing $\rho$ the
error in the finite element solution $u_h$ to problem (\ref{bvp-1})
is bounded by (\ref{lem3.6-1}), i.e.,
\beq
\|(u-u_h)'\|_\Omega \le \frac{C \|r\|_\Omega}{N} .
\label{thm3.1-1}
\eeq
\end{thm}

\vspace{5pt}

As mentioned before, the error bound for a uniform mesh also has the same form
given by (\ref{thm3.1-1}). Although a bound like (\ref{thm3.1-1}) for an equidistributing mesh
is useful in some situations such as in proving the existence of equidistributing meshes
in Lemma \ref{lem3.10}, it does not show any advantage of using an adaptive mesh over
a uniform one.
In the following we shall derive a sharper bound based on the a posteriori error bound (\ref{eta-4}).
The key is to estimate $\alpha_h$, and that is done in a series of lemmas.

\begin{lem}
\label{lem3.7} {\em (Power Inequalities)}
\begin{description}
\item[(i)] Given a real number $0 < q \le 1$, for any  $x,\,y\in \Re$,
 \bey
& |x+y|^q \le |x|^q+|y|^q , &
 \label{lem3.7-1}
\\
& \left | |x|^q - |y|^q \right | \le |x-y|^q .& \label{lem3.7-2}
\eey
\item[(ii)] Given a real number $0 < q \le 1$, for any two functions $v$ and $w$
in a function space equipped with a norm $\| \cdot \|$,
 \bey
& \|v+w\|^q \le \|v\|^q+\|w\|^q , &
 \label{lem3.7-3}
\\
& \left | \|v\|^q - \|w\|^q \right | \le \|v-w\|^q .&
\label{lem3.7-4} \eey
\end{description}
\end{lem}

{\bf Proof.} From the triangle and Jensen's inequalities, (\ref{lem3.7-1}) follows from
\[
 |x+y|^q  \le  \left | |x| + |y| \right |^q \le |x|^q + |y|^q .
\]
Inequality (\ref{lem3.7-2}) is obtained by combining the
inequalities \bey && |x|^q \le |(x-y) + y|^q \le |x-y|^q + |y|^q  ,
\nn \\
&& |y|^q \le |(y-x)+x|^q \le |x-y|^q + |x|^q . \nn \eey

Inequalities (\ref{lem3.7-3}) and (\ref{lem3.7-4}) can be proved
similarly. \proofend

\begin{lem}
\label{lem3.8} For any real number $0< q \le 1$ and any mesh
$\pi_h$ for $\Omega$,
\beq
\| v \|_{L^{2q}(\Omega)}^{2q} \le \sum_i h_i^{1-q} \| v \|_{
K_i}^{2 q} = \sum_i h_i \left <v\right >_i^{2 q}
\le |\Omega|^{1-q} \| v \|_{ \Omega}^{2 q} ,\qquad
\forall v \in L^2(\Omega) . \label{lem3.8-1} \eeq
\end{lem}

{\bf Proof.} The estimates follow from
\[
\sum_i h_i \left <v\right >_i^{2 q} =
\sum_i h_i \left (\frac{1}{h_i} \int_{K_i} |v|^2 d x\right )^q
\ge \sum_i h_i \left (\frac{1}{h_i} \int_{K_i} |v|^{2q} d x\right
) = \| v \|_{L^{2q}(\Omega)}^{2q} ,
\]
\[
\left (\frac{1}{|\Omega|}\sum_i h_i \left (\frac{1}{h_i} \| v
\|_{ K_i}^2 \right )^q\right )^{\frac 1 q} \le
\frac{1}{|\Omega|}\sum_i h_i \left (\frac{1}{h_i} \| v \|_{
K_i}^2 \right ) = \frac{1}{|\Omega|} \| v \|_{ \Omega}^2 .
\]
\proofend

\begin{lem}
\label{lem3.9} For any real number $ 0 < q \le \frac 1 2$ and any mesh
$\pi_h$ for $\Omega$,
\beq
\| v \|_{L^{2q}(\Omega)}^{2q} \le \sum_i h_i^{1-q} \| v \|_{ K_i}^{2
q} \le \| v \|_{L^{2q}(\Omega)}^{2q} + 2 h^{2q} |\Omega|^{1-2q}
\|v'\|_{L^1(\Omega)}^{2q}, \qquad \forall v \in L^2(\Omega),\; v' \in L^1(\Omega)
\label{lem3.9-1}
\eeq
where $h = \max_i h_i$.
\end{lem}

\vspace{5pt}

{\bf Proof.}  The left inequality is a consequence of Lemma
\ref{lem3.8}.

To prove the right inequality, define the element-wise average of
$v$ as
\[
v_{K_i}=\frac1{h_i}\int_{K_i}v dx.
\]
Then, from Lemma \ref{lem3.7}
\bey && \sum_i h_i^{1-q} \| v \|_{ K_i}^{2 q} - \| v
\|_{L^{2q}(\Omega)}^{2q}
\nn \\
& = & \sum_i h_i^{1-q} \|v - v_{K_i} + v_{K_i} \|_{K_i}^{2
q } - \sum_i \int_{K_i} |v|^{2 q} d x
\nn \\
& \le & \sum_i h_i^{1-q} \|v - v_{K_i}\|_{K_i}^{2q} + \sum_i
h_i^{1-q} \| v_{K_i}\|_{K_i}^{2q} - \sum_i \int_{K_i} |v|^{2q}
d x
\nn \\
& = & \sum_i h_i^{1-q} \|v - v_{K_i}\|_{K_i}^{2q} + \sum_i
\int_{K_i} (|v_{K_i}|^{2q} - |v|^{2q}) d x
\nn \\
& \le & \sum_i h_i^{1-q} \|v - v_{K_i}\|_{K_i}^{2q} + \sum_i
\int_{K_i} |v_{K_i} - v|^{2q} d x
\nn \\
& \le & \sum_i h_i^{1-q} \|v - v_{K_i}\|_{K_i}^{2q} + \sum_i
h_i \left (\frac{1}{h_i} \int_{K_i} |v_{K_i} - v|^{2} d x
\right )^q
\nn \\
& \le & 2 \sum_i h_i^{1-q} \|v - v_{K_i}\|_{K_i}^{2q} .
\label{lem3.9-2}
\eey
From the assumption $v'\in L^1(\Omega)$, we
have
\begin{eqnarray}
\|v-v_{K_i}\|^2_{K_i}
&=&\int_{K_i}\left |v(x) - \frac{1}{h_i} \int_{K_i}v(t) d t \right |^2 dx\no\\
&=&\frac1{h_i^2}\int_{K_i}\left |\int_{K_i}(v(x)-v(t))dt\right |^2dx\no\\
&=&\frac1{h_i^2}\int_{K_i}\left |\int_{K_i}\int^x_t v'(s)ds dt\right |^2dx\no\\
&\le & h_i \|v'\|_{L^1(K_i)}^2. \label{lem3.9-3}
\end{eqnarray}
Combining (\ref{lem3.9-3}) with (\ref{lem3.9-2}) and using H\"older's inequality we get
\bey \sum_i
h_i^{1-q} \| v \|_{ K_i}^{2 q} - \| v \|_{L^{2q}(\Omega)}^{2q}
& \le & 2 \sum_i h_i \|v'\|_{L^1(K_i)}^{2q}
\nn \\
& \le & 2 h^{2q}  \sum_i h_i^{1-2q} \|v'\|_{L^1(K_i)}^{2q}
\nn \\
& \le & 2 h^{2q} |\Omega|^{1-2q} \|v'\|_{L^1(\Omega)}^{2q},
\nn
\eey
which gives the right inequality of (\ref{lem3.9-1}).
\proofend

\vspace{15pt}

{\bf Proof of Theorem \ref{thm2.1}.} The bound (\ref{thm2.1-1}) is given by (\ref{eta-4}).

For (\ref{thm2.1-3}), from (\ref{thm2.1-1}) and Lemmas \ref{lem3.3} and \ref{lem3.8} and we have
\bey
\sum_i h_i^{\frac 2 3} \| r - r_h
\|_{K_i}^{\frac 2 3 } & \le &  \|r-r_h\|_\Omega^{\frac 2 3 } 
\nn \\
& \le & C \left [  \| (u-u_h)'\| _\Omega + \| u - u_h \|_\Omega \right ]^{\frac 2 3}
\nn \\
& \le & C \| (u-u_h)'\| _\Omega^{\frac 2 3}
\label{thm2.1-5}
\\
& \le & C N^{- \frac 2 3} \alpha_h^{\frac 1 3} .
\label{thm2.1-6}
\eey
Then from (\ref{alpha-2}), (\ref{thm2.1-6}), and Lemmas \ref{lem3.7} and \ref{lem3.8}
it follows that
\bey \alpha_h^{\frac 1 3} & = &
\sum_i h_i^{\frac 2 3} \| r_h \|_{K_i}^{\frac 2 3 }
\nn \\
& \ge & \sum_i h_i^{\frac 2 3} \| r \|_{K_i}^{\frac 2 3 } - \sum_i
h_i^{\frac 2 3} \| r - r_h \|_{K_i}^{\frac 2 3 }
\nn \\
& \ge & \| r \|_{L^{\frac 2 3}(\Omega)}^{\frac 2 3} - \sum_i
h_i^{\frac 2 3} \| r - r_h \|_{K_i}^{\frac 2 3 } 
\nn \\
& \ge & \| r \|_{L^{\frac 2 3}(\Omega)}^{\frac 2 3} - 
C N^{- \frac 2 3} \alpha_h^{\frac 1 3} ,
\nn \eey
which leads to the left inequality of (\ref{thm2.1-3}) (with $c = C^{\frac 2 3}$).
From (\ref{h-1}), (\ref{thm2.1-6}), Lemmas \ref{lem3.7} and \ref{lem3.9} we have
\bey
\alpha_h^{\frac 1 3} & = & \sum_i h_i^{\frac 2 3} \| r_h\|_{K_i}^{\frac 2 3 }
\nn \\
& \le & \sum_i h_i^{\frac 2 3} \| r \|_{K_i}^{\frac 2 3 } + \sum_i
h_i^{\frac 2 3} \| r - r_h \|_{K_i}^{\frac 2 3 }
\nn \\
& \le & \sum_i h_i^{\frac 2 3} \| r \|_{K_i}^{\frac 2 3 } + 
C N^{- \frac 2 3} \alpha_h^{\frac 1 3}
\nn \\
& \le &  \| r \|_{L^{\frac 2 3}(\Omega)}^{\frac 2 3} + 
C N^{- \frac 2 3} \| r' \|_{L^1(\Omega)} ^{\frac 2 3} +
C N^{- \frac 2 3} \alpha_h^{\frac 1 3} ,
\nn
\eey
which yields the right inequality of (\ref{thm2.1-3}).

We now prove (\ref{thm2.1-2}). In this situation, $r \in L^2(\Omega)$.
From the above derivation we can see that, for $N > C^{\frac 3 2}$,
\beq
\left (1+ C N^{- \frac 2 3} \right )^{-1} \|r \|_{L^{\frac 2 3}(\Omega)}^{\frac 2 3} \le
\alpha_h^{\frac 1 3}  \le \left (1- C N^{- \frac 2 3} \right )^{-1} \sum_i h_i^{\frac 2 3} \| r
\|_{K_i}^{\frac 2 3 } .
\label{thm2.1-7}
\eeq
Since functions having derivatives in $L^1(\Omega)$ are
dense in $L^2(\Omega)$, given any $\e>0$  there
exists a function $\tilde{r}$ such that
 \beq
\tilde{r}' \in L^1(\Omega)\quad \mbox{ and }\quad
\|r-\tilde{r}\|_{\Omega}\le \e. \label{thm2.1-8} \eeq
Then, from Lemmas \ref{lem3.7}, \ref{lem3.8}, and \ref{lem3.9}) we have
\bey &&
\sum_i h_i^{\frac{2}{3}}\|r\|_{K_i}^{\frac 2 3} - \int_{\Omega}
|r|^{\frac{2}{3}} d x
\nn \\
& = & \sum_i h_i^{\frac{2}{3}}\|r-\tilde{r}+\tilde{r}\|_{K_i}^{\frac
2 3} - \int_{\Omega} |\tilde{r}|^{\frac{2}{3}} d x + \left (
\int_{\Omega} |\tilde{r}|^{\frac{2}{3}} d x - \int_{\Omega}
|r|^{\frac{2}{3}} d x\right )
\nn \\
& \le & \sum_i h_i^{\frac{2}{3}}\left ( \|r-\tilde{r}\|_{K_i}^{\frac
2 3} +\| \tilde{r}\|_{K_i}^{\frac 2 3}\right ) - \int_{\Omega}
|\tilde{r}|^{\frac{2}{3}} d x + \int_{\Omega}
|\tilde{r}-r|^{\frac{2}{3}} d x
\nn \\
& = & \left (\sum_i h_i^{\frac{2}{3}}\|\tilde{r}\|_{K_i}^{\frac 2 3}
- \int_{\Omega} |\tilde{r}|^{\frac{2}{3}} d x\right ) + \sum_i h_i
\left (\frac{1}{h_i} \|r-\tilde{r}\|_{K_i}^2\right )^{\frac 1 3} +
\int_{\Omega} |\tilde{r}-r|^{\frac{2}{3}} d x
\nn \\
& \le & C N^{-\frac 2 3} \| \tilde{r}' \|_{L^1(\Omega)}^{\frac 2 3} + 2 \|r - \tilde{r}\|_{\Omega}^{\frac 2 3}
\nn \\
& \le & C N^{-\frac 2 3} \| \tilde{r}' \|_{L^1(\Omega)}^{\frac 2 3} + 2  \e^{\frac 2 3} .
\nn
\eey
Inserting this into (\ref{thm2.1-7}) gives
\beq
\left (1+ C N^{-\frac 2 3} \right )^{-1} \|r
\|_{L^{\frac 2 3}(\Omega)}^{\frac 2 3} \le \alpha_h^{\frac 1 3}  \le
\left (1- C N^{-\frac 2 3} \right )^{-1}
\left [ \|r \|_{L^{\frac 2 3}(\Omega)}^{\frac 2 3} + C N^{-\frac 2 3} \| \tilde{r}' \|_{L^1(\Omega)}^{\frac2 3}
 + 2 \e^{\frac 2 3} \right ].
\label{thm2.1-9}
\eeq
Taking limit as $N \to \infty$ in the above inequality yields
\[
\|r \|_{L^{\frac 2 3}(\Omega)}^{\frac 2 3} \le \lim_{N \to \infty}
\alpha_h^{\frac 1 3} \le \|r \|_{L^{\frac 2 3}(\Omega)}^{\frac 2 3}
+ 2 \e^{\frac 2 3} .
\]
Finally, taking limit as $ \epsilon \to 0$ in the above inequality
gives (\ref{thm2.1-2}).
\proofend

\vspace{10pt}

It is remarked that we can also use the a priori error bound (\ref{lem3.6-1})
to estimate $\|(u-u_h)'\|_\Omega$ in (\ref{thm2.1-5}).
%In this situation, the mesh does not have to be equidistributing but
%$u\in H^2(\Omega)$ is required (cf. Lemmas \ref{lem3.4} and \ref{lem3.6}).
For convenience, we list the result in the following lemma without giving the detail
of the proof.

\begin{lem}
\label{lem3.10}
Assume that the solution to problem (\ref{bvp-1}) satisfies $u \in H^2(\Omega)$
and the mesh has the property (\ref{h-1}).  If $r$ is only $L^2$ integrable,
$\alpha_h$ has the property (\ref{thm2.1-2}). If further $r' \in L^1(\Omega)$, then
$\alpha_h$ defined in (\ref{alpha-2}) is bounded by
\beq
\|r\|_{L^{\frac 2 3}(\Omega)}^{\frac 2 3} - \left (\frac{C \|r\|_\Omega}{N}\right )^{\frac 2 3}
\;\le\; \alpha_h^{\frac 1 3} \;\le\;
\| r \|_{L^{\frac 2 3}(\Omega)}^{\frac 2 3}
+ \left (\frac{C}{N} \| r' \|_{L^1(\Omega)} \right )^{\frac 2 3} + \left (\frac{C \|r\|_\Omega}{N}\right )^{\frac 2 3} .
\label{lem3.10-1}
\eeq
\end{lem}

\vspace{10pt}

{\bf Proof of Theorem \ref{thm2.2}.} When the adaptation function
$\rho$ and the intensity parameter $\alpha_h$ are chosen as in
(\ref{rho-2}) and (\ref{alpha-2}), a quasi-equidistributing mesh
satisfying (\ref{eq-3}) has the property
\beq
h_i \le \frac{2\kappa}{N},\qquad i = 1, ...,N.
\label{h-3}
\eeq
Moreover, from (\ref{eq-3}), (\ref{eta-1}), and (\ref{sigma-2})  we have
\bey
\tilde{\eta}_h^2 & = &
\alpha_h \sum_{i} \left (h_i \rho_{i} \right)^3
\nn \\
&\le & \alpha_h \sum_{i} \left
(\frac{\kappa \sigma_h}{N} \right )^3 
\nn \\
& = & \frac{\alpha_h (\kappa \sigma_h)^3}{N^2}
\nn \\
& \le & \frac{8 \alpha_h \kappa^3 }{N^2} .
\nn
\eey
The remaining of the theorem can be proven similarly as for Theorem \ref{thm2.1}. \proofend

% section 4
\section{An iterative algorithm for computing equidistributing meshes and numerical examples}
\label{SEC:numerical-results}

We start with describing an iterative algorithm for computing
equidistributing meshes. Recall that the finite element equation (\ref{fem-1}), the
equidistribution relation (\ref{eq-2}), and the boundary conditions
$x_0 = 0$ and $x_N=1$ form a nonlinear algebraic system for the
physical solution $u_h$ and the mesh $\pi_h$. This system is
typically solved iteratively; see Fig. \ref{f1}. An algorithm of this type
is given in the following. Starting from an initial mesh $\pi_h^{(0)}$,
it produces a sequence of meshes and solutions, $\{ \pi_h^{(k)}, u_h^{(k)}\} $.

\vspace{10pt}

\noindent {\bf Algorithm for computing
equidistributing meshes.} Given an integer $N > 0$ and an initial
mesh $\pi_h^{(0)}$, for $k=0, 1, ...$ do
\begin{itemize}
\item[(i)] Solution of the boundary value problem using mesh $\pi_h^{(k)}$. This step is to find $u_h^{(k)}
\in V_h^{(k)}$ such that \beq B(u_h^{(k)},v_h)=(f,v_h), \quad
\forall\, v_h\in V_h^{(k)} . \label{fem-2} \eeq

\item[(ii)]  Mesh generation. This step is to compute the new equidistributing mesh
using the equidistribution relation (\ref{eq-2}), i.e., \beq
\int_{0}^{x_i^{(k+1)}} \rho^{(k)}(x) d x =
\frac{i}{N}\sigma_h^{(k)},\quad i = 1, ..., N-1 \label{eq-4} \eeq
where \bey && \rho^{(k)}(x) =  \rho_i^{(k)} \equiv \left
(1+\frac{1}{\alpha_h^{(k)} } \left <r_h^{(k)}\right >_{K_i^{(k)}}^2\right )^{\frac 1 3},\quad \forall x \in
K_i^{(k)}, \quad i = 1, ..., N \label{rho-3}
\\
&&  \alpha_h^{(k)} = \left [  \sum_i h_i^{(k)} \left <r_h^{(k)}\right >_{K_i^{(k)}}^{\frac
2 3} \right ]^3 , \label{alpha-3}
\\
&& \sigma_h^{(k)} = \sum_{i=1}^N
h_i^{(k)} \rho_i^{(k)} . \label{sigma-3} \eey
\end{itemize}

\vspace{10pt}

Note that the left-hand side of (\ref{eq-4}) is a
monotone, piecewise linear function of $x_i^{(k+1)}$ and
an explicit formula can be found as
\[
x_{i}^{(k+1)} = x_{j-1}^{(k)} + \frac{\left( \frac{i}{N}
\sigma_h^{(k)} - \sum_{l=1}^{j-1} h_l^{(k)} \rho_l^{(k)}
\right)}{\rho_j^{(k)}},
\]
where $j$ is the index satisfying
\[
\sum_{l=1}^{j-1} h_l^{(k)} \rho_l^{(k)} <  \frac{i}{N} \sigma_h^{(k)} \leq
\sum_{l=1}^{j} h_l^{(k)} \rho_l^{(k)} .
\]
Moreover, Steps (i) and (ii) define a map
\beq
G_N: \Re^{N-1} \to \Re^{N-1}:\quad 
X^{(k+1)} = G_N X^{(k)},
\label{GN}
\eeq
where $X^{(k)}$ and $X^{(k+1)}$ are the $(N-1)$-component vectors corresponding to the meshes
$\pi_h^{(k)}$ and $\pi_h^{(k+1)}$, respectively; see (\ref{SN}). It is not difficult to see that a fixed point of this map satisfies
(\ref{eq-1}) and thus is an equidistributing mesh.
Furthermore, the computation can be stopped when
\beq \|\pi_h^{(k+1)}-\pi_h^{(k)}\|_\infty \equiv \max_i |x_i^{(k+1)} -
x_i^{(k)}| \le \epsilon \label{stop-1} \eeq or \beq Q_{eq,i}^{(k)}
\equiv \frac{N \rho_i^{(k)} h_i^{(k)}}{\sigma_h^{(k)}} \le
\kappa, \quad i = 1, ..., N \label{stop-2} \eeq where $\epsilon
> 0$ is a prescribed tolerance, $\kappa $ is a number chosen to be
close to and greater than one, and $Q_{eq,i}^{(k)}$ is the so-called
quality measure of equidistribution \cite{Hua03}. The second
stopping criterion needs some explanation.  It is not difficult to see
that $Q_{eq,i}^{(k)}$ has the properties
\beq
\frac{1}{N} \sum_i Q_{eq,i}^{(k)} = 1,\qquad \max_i Q_{eq,i}^{(k)}
\ge 1.
\label{Qeq-1}
\eeq
In addition,  $\max_i Q_{eq,i}^{(k)} = 1$ if and only if the mesh is
an equidistributing mesh satisfying (\ref{eq-2}). Thus, if
the mesh sequence $\pi_h^{(k)}$ converges to an equidistributing mesh
we will have $ \max_i Q_{eq,i}^{(k)} \to 1$ as $k\to \infty$; and vice versa.
This implies that (\ref{stop-2}) is an effective stopping criterion.

It is interesting to point out that
$\max_i Q_{eq,i}^{(k)}$ actually measures how closely the equidistribution relation
(\ref{eq-2}) is satisfied by the mesh; see \cite{Hua03} for detailed discussion.
Moreover, by the definition (\ref{eq-3}) one can see that 
any mesh satisfying (\ref{stop-2}) is a quasi-equidistributing mesh.
Finally, from (\ref{stop-2})  and (\ref{Qeq-1}) we have
\[
- N (\kappa-1) + \kappa \le Q_{eq,i}^{(k)} \le \kappa, \quad i = 1,
..., N
\]
where $ - N (\kappa-1) + \kappa > 0 $ when $\kappa$ is sufficiently close to one.

\vspace{5pt}

We now present some numerical results to demonstrate the convergence of the algorithm.

\vspace{5pt}

{\bf Example 3.1.} This example is a reaction-diffusion equation
\beq -\epsilon u'' + u = -2\epsilon-x (1-x)-1 \label{exam-1-1} \eeq
subject to the boundary condition (\ref{bc-1}). The exact solution is
given by
\beq u =
\frac{1}{1-e^{-\frac{2}{\sqrt{\epsilon}}}} \left
(e^{-\frac{1-x}{\sqrt{\epsilon}}}-e^{-\frac{1+x}{\sqrt{\epsilon}}} +
e^{-\frac{x}{\sqrt{\epsilon}}}
        - e^{-\frac{2-x}{\sqrt{\epsilon}}} \right ) - x(1-x)-1 .
\eeq It exhibits boundary layers at both ends of interval $[0,1]$  when $\epsilon$
is small. The parameter is taken as $\epsilon = 10^{-5}$.

\begin{figure}[thb]
\centering
\hbox{
\begin{minipage}[t]{3in}
\includegraphics[width=3in]{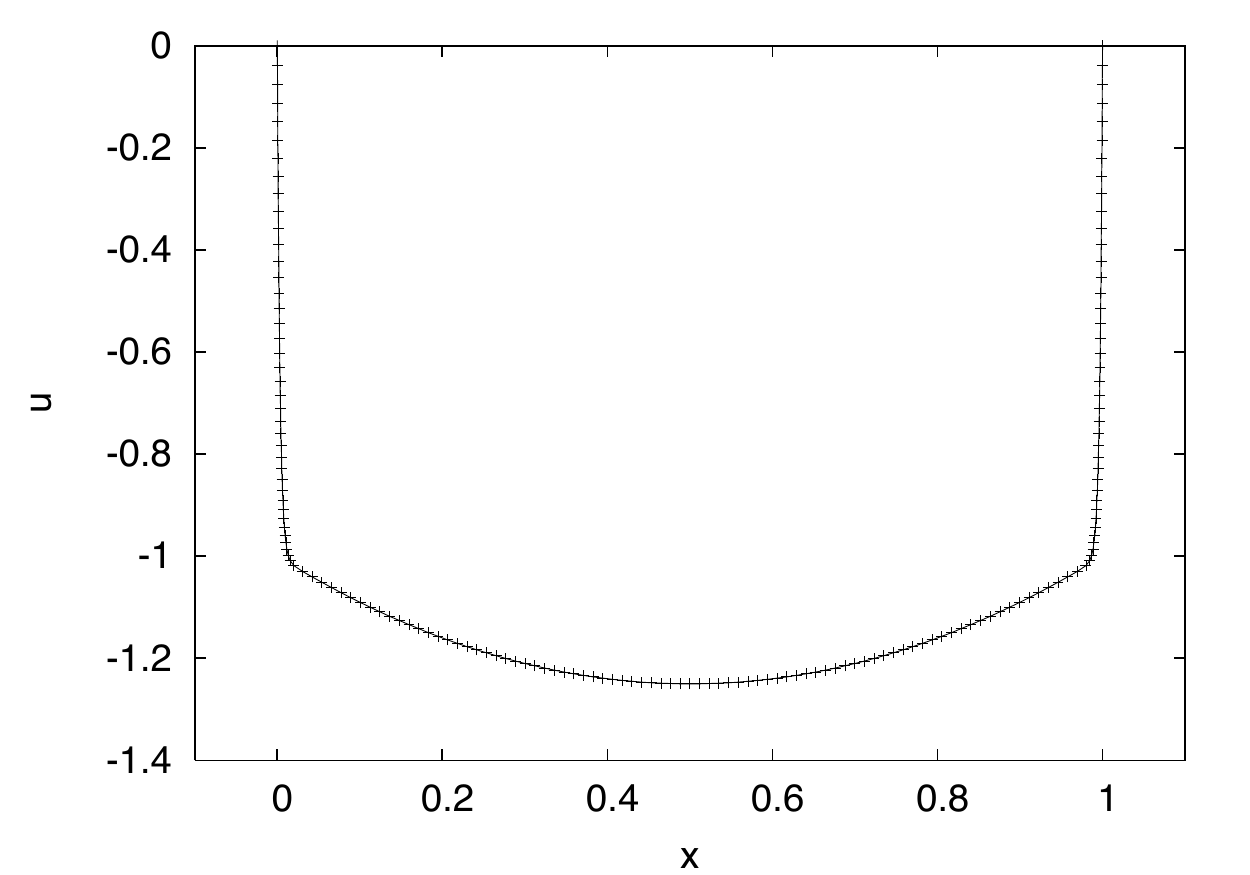}
\end{minipage}
\hspace{10mm}
\begin{minipage}[t]{3in}
\includegraphics[width=3in]{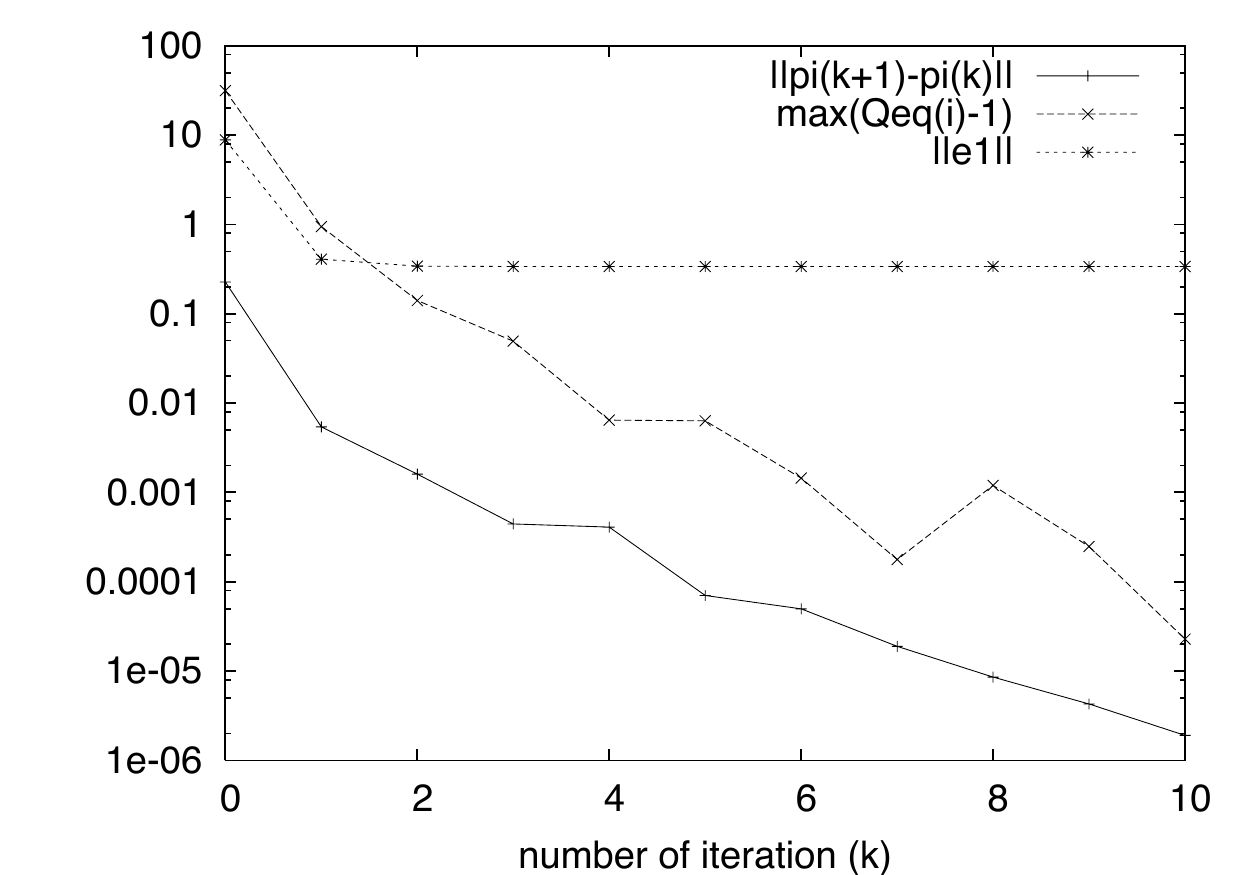}
\end{minipage}
}
\caption{Example 3.1. (a) An adaptive mesh of $N=161$ points is plotted on the curve of the computed solution.
(b) The difference between consecutive meshes ($\|\pi_h^{(k+1)}-\pi_h^{(k)}\|_\infty$), the equidistribution quality
measure ($\max_i (Q_{eq,i}-1)$), and the solution error ($\|(u_h-u)'\|_\Omega$) are plotted against the number
of iteration, $k$.}
\label{fex1-1}
\end{figure}

A typical adaptive mesh and the corresponding computed solution are
shown in Fig. \ref{fex1-1}(a). In Fig. \ref{fex1-1}(b),
$\|\pi_h^{(k+1)}-\pi_h^{(k)}\|_\infty$,  $\max_i (Q_{eq,i}-1)$, and
$\|(u_h-u)'\|_\Omega$ are plotted as functions of the number of
iteration. It can be seen that
 both (\ref{stop-1}) and (\ref{stop-2}) are effective stopping criteria and
$\|\pi_h^{(k+1)}-\pi_h^{(k)}\|_\infty$ and $\max_i (Q_{eq,i}-1)$ converge in
a similar manner. Moreover, the solution error quickly reaches its
lowest level (in one or two iterations for the current case).

The number of iterations required to reach the stopping criterion
$\max_i Q_{eq,i} \le 1.01$ and the solution error and the modified a
posteriori estimator on the final mesh of each run are listed in
Table \ref{tex1-1}. The results show that the underlying iterative
algorithm may fail for small $N$ but is convergent for sufficiently
large $N$. Moreover, the algorithm converges faster for larger $N$. These
results are consistent with the observations made in Pryce
\cite{Pry89} and Xu et al. \cite{XHRW08} for the convergence of de Boor's
algorithm for generating
equidistributing meshes for a given analytical function. It can also
be seen that $\|(u-u_h)'\|_\Omega$ is smaller than the error estimator $\tilde{\eta}$ and both
$\|(u-u_h)'\|_\Omega$ and $\tilde{\eta}$ converge in the same order $O(1/N)$
as $N\to \infty$. These results conform the theoretical predictions in Theorem
\ref{thm2.1} and \ref{thm2.2}.

% table tex1-1
\begin{table}[htb]
\begin{center}
\caption{Example 3.1. $Iter$ is the number of iterations required to reach the stopping criterion
$\max_i Q_{eq,i} \le 1.01$ or the maximum allowed number (1000 is used in the computation).
$\|(u_h-u)'\|_\Omega$ is the error obtained for the final mesh of each computation.}
\vspace{5pt}
\begin{tabular}{c|llllll}\hline \hline
$N$ & 21 & 41 & 81 & 161 & 321 & 641 \\ \hline
$Iter$ & 1000 & 1000 & 39 & 4 & 3 & 2  \\ \hline
$\|(u_h-u)'\|_\Omega$ & 3.07 & 1.41 & 6.86e-1 & 3.39e-1 & 1.69e-1 & 8.46e-2 \\ \hline
$\tilde{\eta}$ & 8.10 & 3.65 & 1.72 & 8.35e-1 & 4.15e-1 & 2.07e-1  \\
\hline \hline
\end{tabular}
\label{tex1-1}
\end{center}
\end{table}

% fex1-1

\vspace{10pt}

\begin{figure}[thb]
\centering
\hbox{
\begin{minipage}[t]{3in}
\includegraphics[width=3in]{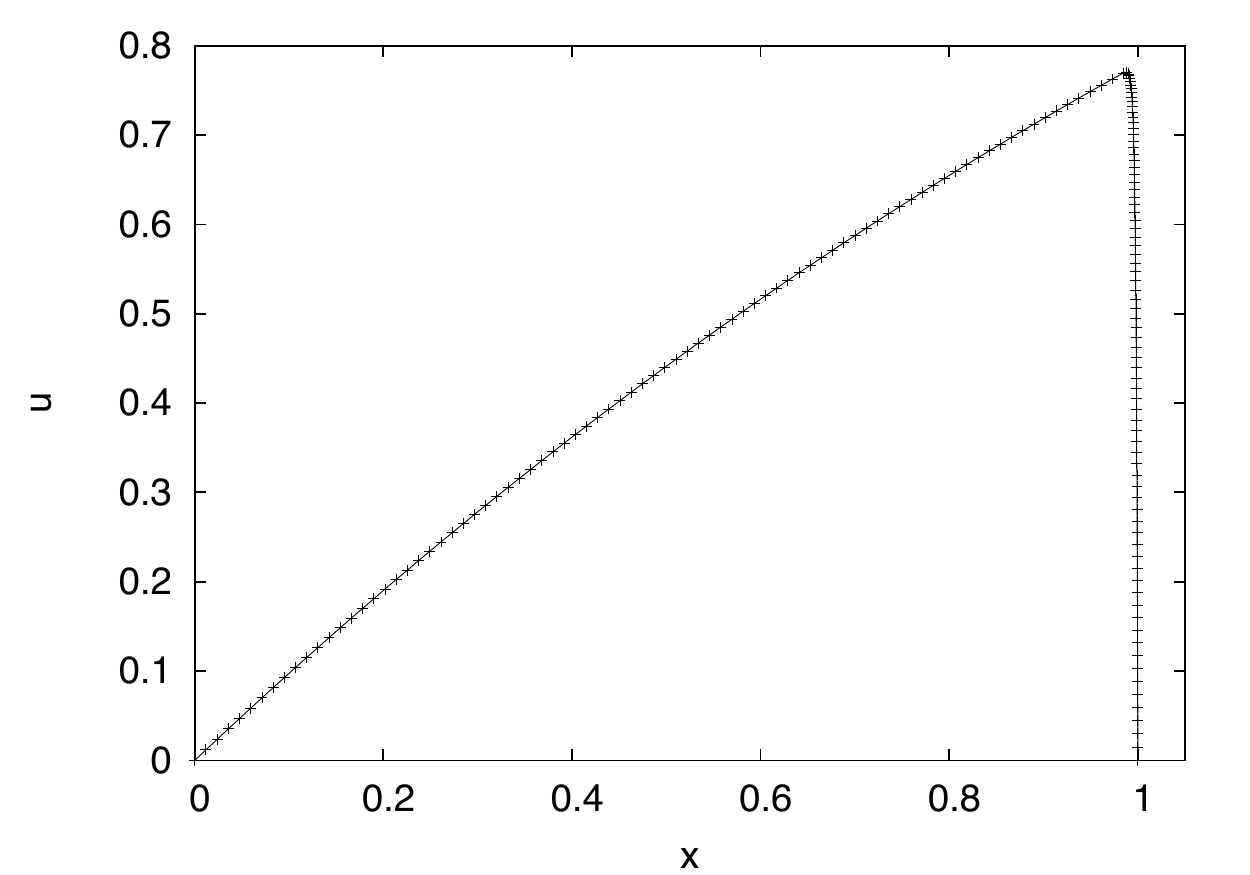}
\end{minipage}
\hspace{10mm}
\begin{minipage}[t]{3in}
\includegraphics[width=3in]{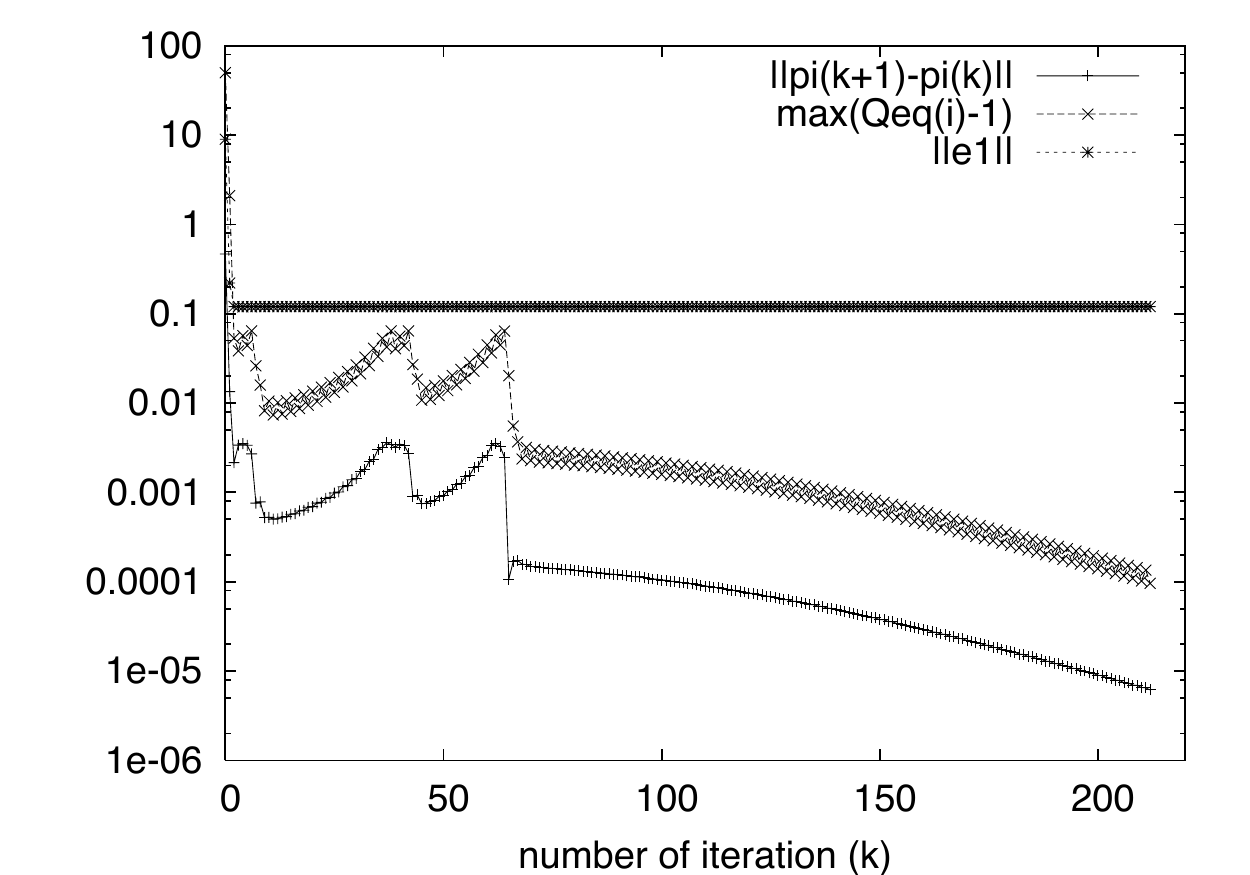}
\end{minipage}
}
\caption{Example 3.2. (a) An adaptive mesh of $N=161$ points is plotted on the curve of the computed solution.
(b) The difference between consecutive meshes ($\|\pi_h^{(k+1)}-\pi_h^{(k)}\|_\infty$), the equidistribution quality
measure ($\max_i (Q_{eq,i}-1)$), and the solution error ($\|(u_h-u)'\|_\Omega$) are plotted against the number
of iteration, $k$.}
\label{fex2-1}
\end{figure}

{\bf Example 3.2.} Our second example is a convection-dominated differential equation
\beq
-\epsilon u'' + (1-\frac{ 1}{ 2}\epsilon) u' + \frac{1}{4}\left (1-\frac{1}{4}\epsilon\right ) u =  e^{-\frac{x}{4}} ,
\label{exam-2-1}
\eeq
where $\epsilon = 2\times 10^{-3}$. The exact solution is given by
\beq
u = e^{-\frac{x}{4}} \left ( x - \frac{e^{-\frac{1-x}{\epsilon}} - e^{-\frac{1}{\epsilon}}}{1-e^{-\frac{1}{\epsilon}}}  \right ),
\eeq
which has the boundary layer at $x=1$ when $\epsilon$ is small. The numerical results are
showed in Fig. \ref{fex2-1} and Table \ref{tex2-1}. These results confirm the observations made
from the previous example. Particularly, the algorithm converges for sufficiently large $N$
and faster for larger $N$. Moreover, the $H^1$ semi-norm of
the error converges in the first order $O(1/N)$ as $N\to \infty$.

% table tex2-1
\begin{table}[htb]
\begin{center}
\caption{Example 3.2. $Iter$ is the number of iterations required to reach the stopping criterion
$\max_i Q_{eq,i} \le 1.01$ or the maximum allowed number (1000 is used in the computation).
$\|(u_h-u)'\|_\Omega$ is the error obtained for the final mesh of each computation.}
\vspace{5pt}
\begin{tabular}{c|llllll}\hline \hline
$N$ & 21 & 41 & 81 & 161 & 321 & 641 \\ \hline
$Iter$ & 1000 & 327 & 83 & 9 & 5 & 3  \\ \hline
$\|(u_h-u)'\|_\Omega$ & 1.20 & 5.07e-1 & 2.54e-1 & 1.20e-1 & 5.95e-2 & 2.96e-2 \\ \hline
$\tilde{\eta}$ & 8.08 & 1.30 & 6.44e-1 & 2.97e-1 & 1.46e-1 & 7.26e-2 \\
\hline \hline
\end{tabular}
\label{tex2-1}
\end{center}
\end{table}

% fex2-1

\vspace{10pt}

\begin{figure}[thb]
\centering
\includegraphics[width=3in]{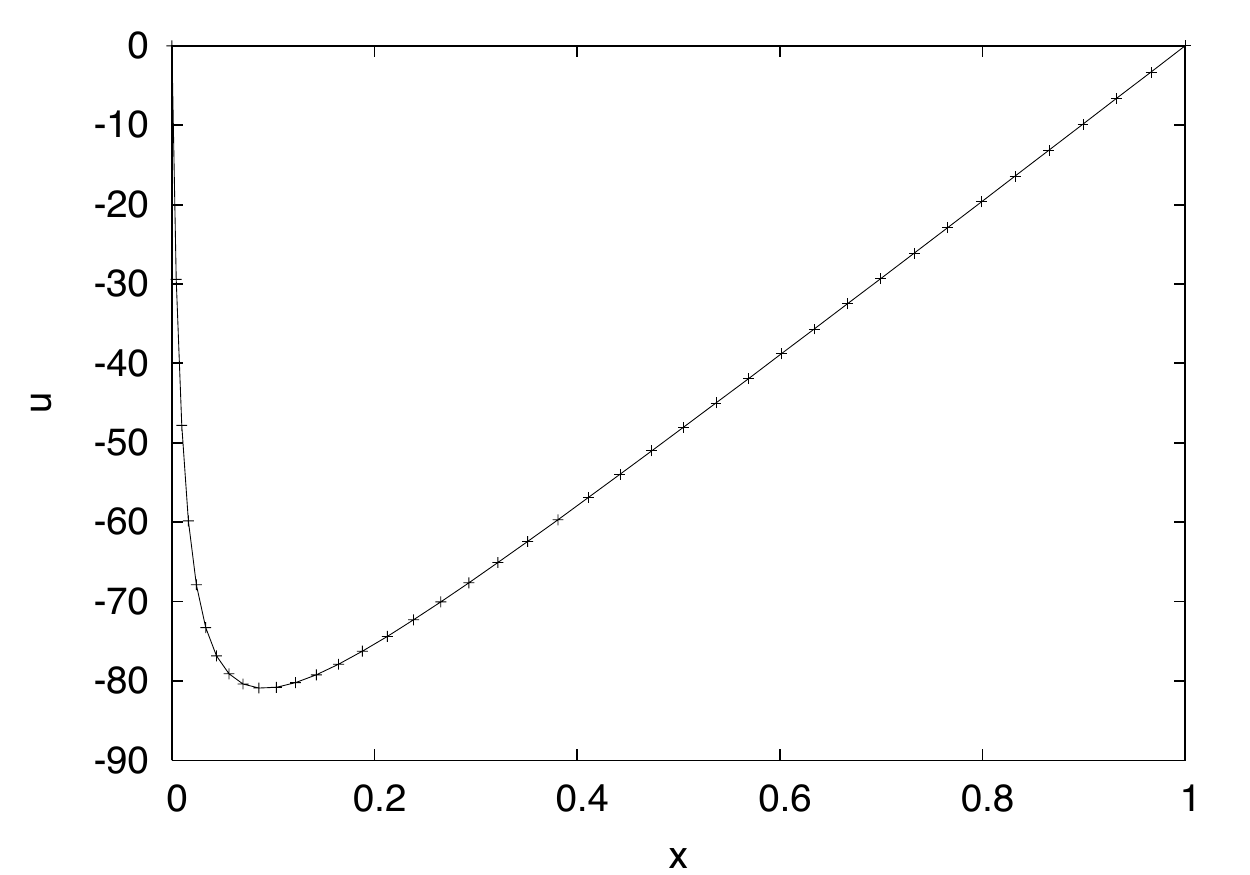}
\caption{Example 3.3. An adaptive mesh of $N=41$ points is plotted on the curve of the computed solution.}
\label{fex3-1}
\end{figure}

{\bf Example 3.3.} This example has been used by Babu\u{s}ka and Rheinboldt \cite{BR79b}.
It takes the form
\beq
- ((x+\alpha)^p u')' + (x+\alpha)^q u = f,
\eeq
where $f$ is chosen such that the exact solution of the boundary value problem (with boundary condition
(\ref{bc-1})) is
\beq
u = (x+\alpha)^r - \left ( \alpha^r (1-x)+(1+\alpha)^r x \right ).
\eeq
In our computation, the parameters are taken as $p = 2$, $q = 1$, $r = -1$, and $\alpha = 1/100$.
The numerical results are shown in Fig. \ref{fex3-1} and Table \ref{tex3-1}. Once again, these results
confirm the observations made from the previous examples.

% table tex3-1
\begin{table}[htb]
\begin{center}
\caption{Example 3.3. $Iter$ is the number of iterations required to reach the stopping criterion
$\max_i Q_{eq,i} \le 1.01$ or the maximum allowed number (1000 is used in the computation).
$\|(u_h-u)'\|_\Omega$ is the error obtained for the final mesh of each computation.}
\vspace{5pt}
\begin{tabular}{c|llllll}\hline \hline
$N$ & 21 & 41 & 81 & 161 & 321 & 641 \\ \hline
$Iter$ & 4 & 3 & 3 & 2 & 2 & 2  \\ \hline
$\|(u_h-u)'\|_\Omega$ & 2.15e2 & 1.13e2 & 5.73e1 & 2.88e1 & 1.44e1 & 7.20 \\ \hline
$\tilde{\eta}$ & 3.22e3 & 1.51e3 & 7.41e2 & 3.69e2 & 1.84e2 & 9.21e1 \\
\hline \hline
\end{tabular}
\label{tex3-1}
\end{center}
\end{table}

% fex3-1

% section 5
\section{Continuous dependence of finite element solution on mesh}
\label{SEC:continuity}

This section is devoted to the proof of Theorem \ref{thm2.3} for the continuous
dependence of the linear finite element solution on mesh. To this end, we need to establish
some error bounds in the $L^\infty$ norm, which are also needed in the next section
in obtaining the upper and lower bounds for the adaptation function. 
It is worth pointing out that in this section we do not require that the mesh be necessarily
an equidistributing mesh. Instead, all results hold for any mesh in $S_N$.
Especially, Lemmas \ref{lem3.6} and \ref{lem3.10} are true for any mesh having property (\ref{h-1}).

We shall use two different meshes, $\pi_h$ and $\pi_{\tilde{h}}$ or
$X$ and $\tilde{X}$, in this and next sections. To distinguish the dependence we shall denote
any quantity or function (say $v$) associated with mesh $\tilde{X}$ by $\tilde{v}$.
Moreover, in these two sections constants are considered  as numbers that may further depend
on the solution $u$ and the right-hand side function $f$ (but not on the mesh).

We start with establishing two inequalities in Lemma \ref{lem5.1} and
error bounds in the $L^\infty$ norm  in Lemma \ref{lem5.2}. 

\begin{lem}
\label{lem5.1}
\bey
|v(x)|^2&\le&\int^{b_1}_{a_1} |v|~|v'|dt,\quad \forall x \in [a_1, b_1],\quad \forall v\in H^1_0(a_1,b_1)
\label{lem5.1-1}
\\
|v(x)|^2&\le&2\int^{b_1}_{a_1}|v|~|v'|dt+\frac1{b_1-a_1}\int^{b_1}_{a_1}|v|^2dt ,
\quad \forall x \in [a_1, b_1], \quad \forall v\in H^1(a_1,b_1) .
\label{lem5.1-2}
\eey
\end{lem}

{\bf Proof.} These results can readily be proven using integration by parts.
%For $v\in H^1_0(a_1,b_1)$ we have
%\bey
%&& |v(x)|^2=2\int^x_{a_1}v~v'dt \le 2 \int_{a_1}^x |v| \cdot |v'| d t,
%\nn \\
%&& |v(x)|^2=-2\int^{b_1}_x v~v'dt \le 2 \int_x^{b_1} |v| \cdot |v'| d t.
%\nn
%\eey
%Summing these inequalities yields (\ref{lem5.1-1}).

%Moreover, for $v\in H^1(a_1,b_1)$, using integration by parts we have
%\bey &&
%\int^{x}_{a_1}|v|^2dt=(x-a_1)|v(x)|^2-2\int^x_{a_1}(t-a_1)v~v'dt \ge
%(x-a_1)|v(x)|^2-2(b_1-a_1)\int^x_{a_1} |v|\cdot |v'| dt,
%\nn \\
%&& \int^{b_1}_{x}|v|^2dt=(b_1-x)|v(x)|^2-2\int^{b_1}_x(t-b_1)v~v' dt
%\ge (b_1-x)|v(x)|^2-2(b_1-a_1)\int^{b_1}_x |v|\cdot |v'| dt.
%\no
%\eey
%Summing these two inequalities gives (\ref{lem5.1-2}).
\proofend

\begin{lem}
\label{lem5.2}
Assume that $X \in S_N$ and $u \in H^2(\Omega)$. Then the finite element error can be bounded
in $L^\infty$ norm as
\bey
\|e_h\|_{L^\infty(\Omega)} & \le & \frac{C \|r \|_\Omega}{N},
\label{lem5.2-1}
\\
\|e_h'\|_{L^\infty(\Omega)} & \le & \frac{\sqrt{\rho_0+1} C \|r \|_\Omega}{\sqrt{N}} .
\label{lem5.2-2}
\eey
\end{lem}

\vspace{5pt} 

{\bf Remark.} The dependence on constant $\rho_0$ is spelled out explicitly in (\ref{lem5.2-2}).
This is needed for the definition of $\rho_0$; see the proof of Lemma \ref{lem6.0}.

\vspace{5pt}

{\bf Proof.} From Poincare's inequality and Lemma \ref{lem3.6} 
we know that the error can be bounded in $L^2$ norm as
\beq
\|e_h\|_\Omega \le \frac{C \| r \|_\Omega}{N},\quad
\|e_h'\|_\Omega \le \frac{C \| r \|_\Omega}{ N} .
\label{lem5.2-3}
\eeq
Then from Lemma \ref{lem5.1} and Schwarz' inequality we get
\[
\| e_h \|_{L^\infty (\Omega)}^2 \le \|e_h\|_\Omega \|e_h'\|_\Omega
\le \frac{C \| r \|_\Omega^2}{N^2} ,
\]
which leads to (\ref{lem5.2-1}).

To prove (\ref{lem5.2-2}), taking $v=e_h'$ and $(a_1, b_1) = K_i$ ($i = 1, ..., N$) in (\ref{lem5.1-2}) we have
\beq
\|e'_h\|^2_{L^\infty(K_i)}\le\frac{1}{h_i}\|e'_h\|^2_{K_i}+2\|e'_h\|_{K_i}\|u''\|_{K_i}.
\label{lem5.2-4}
\eeq
Noticing that $1/(\rho_0 N) \le h_i$, we get from (\ref{lem5.2-3}) and (\ref{lem5.2-4}) that
\bey
\|e_h'\|_{L^\infty(K_i)}^2 & \le & \rho_0 N
\|e'_h\|^2_{\Omega}+2\|e'_h\|_{\Omega}\|u''\|_{\Omega}
\nn \\
& \le & \rho_0 N \|e'_h\|^2_{\Omega}+C \|e'_h\|_{\Omega}\|r\|_{\Omega}
\nn \\
& \le &   \frac{(\rho_0 + 1) C \| r \|_\Omega^2}{ N} ,
\nn
\eey
which gives (\ref{lem5.2-2}).
\proofend

\vspace{10pt}

Note that the convergence order for $\|e_h\|_{L^\infty(\Omega)}$ is not optimal in the above lemma.
The optimal can be obtained by making use of the Nitsche trick.

\begin{lem}
\label{lem5.4} Assume that $X \in S_N$ and $u \in H^2(\Omega)$. Then
the finite element error $e_h$ is bounded by
\bey
&& \|e_h\|_\Omega \le \frac{C}{N^2},\qquad \|e'_h\|_\Omega \le \frac{C}{N},
\label{lem5.4-1}
\\
&& \|e_h\|_{L^\infty(\Omega)}\le \frac{C}{N^{\frac 3 2}},\qquad \|e'_h\|_{L^\infty(\Omega)} \le \frac{C}{\sqrt{N}} ,
\label{lem5.4-2}
\\
&& \|u_h\|_{L^\infty(\Omega)}\le C,\qquad \|u_h'\|_{L^\infty(\Omega)} \le C ,
\label{lem5.4-3}
\eey
where the generic constant $C$ may further depend on the solution $u$ and the right-hand side function $f$.
\end{lem}

{\bf Proof.} The second inequality of (\ref{lem5.4-1}) is a consequence of Lemma \ref{lem3.6}.
The first inequality can readily be proven by making use of the Nitsche trick. The proof for (\ref{lem5.4-2})
is similar to that for Lemma \ref{lem5.2} but makes use of (\ref{lem5.4-1}).
 The inequalities in (\ref{lem5.4-3}) follow from (\ref{lem5.4-2}),
the triangle inequality, and the boundedness of $u$ and $u'$ in $L^\infty$ norm.
\proofend

\vspace{10pt}

We now consider the continuous dependence of the finite element solution $u_h$ on mesh.

\begin{lem}
\label{lem5.5}
Assume that $X,\; \tilde{X} \in S_N$, $f \in L^\infty(\Omega)$, and $u \in H^2(\Omega)$.
Then the finite element solutions $u_h$ related to mesh $X$ and $u_{\tilde{h}}$
related to $\tilde{X}$ satisfy
\bey
&&  \|(u_h-\hat{u}_h)'\|^2_\Omega+\|u_h-\hat{u}_{h}\|^2_\Omega+
\sum^N_{i=1}\frac{1}{h_i}|u_h(x_i)-u_{\tilde{h}}(\tilde{x}_i)-(u_h(x_{i-1})-u_{\tilde{h}}(\tilde{x}_{i-1}))|^2
 \no\\
&& \qquad \hfill{} \qquad  + \; \sum^{N}_{i=1} h_i|u_h(x_i)-u_{\tilde{h}}(\tilde{x}_i)|^2  \le
C N^2 \; \|X-\tilde{X}\|_\infty^2,
\label{lem5.5-1}
\eey
where
\[
\hat{u}_h(x) = \sum_{i} u_{\tilde{h}}(\tilde{x}_i) \phi_i(x).
\]
\end{lem}

{\bf Proof.} We notice that the finite element solutions can be expressed as
\[
u_h=\sum\limits_{i=1}^{N-1}u_h(x_i)\phi_i(x),\quad 
u_{\tilde{h}}=\sum\limits_{i=1}^{N-1}u_{\tilde{h}}(\tilde{x}_i)\tilde{\phi}_i(x).
\]
Moreover, (\ref{fem-1}) can be rewritten into matrix form as
\beq
A U = F,  \quad \tilde{A} \tilde{U} = \tilde{F},
\label{lem5.5-2}
\eeq
where
\bey
&& A=(a_{ij}), ~\tilde{A}=(\tilde{a}_{ij}),\qquad a_{ij}=B(\phi_j,\phi_i),
~~\tilde{a}_{ij}=B(\tilde{\phi}_j,\tilde{\phi}_i),
\nn \\
&& U=(u_h(x_1),\cdots, u_h(x_{N-1}))^{\top},
\qquad \tilde{U}=(u_{\tilde{h}}(\tilde{x}_1),\cdots,
u_{\tilde{h}}(\tilde{x}_{N-1}))^{\top},
\nn \\
&& F=(F_1,\cdots, F_{N-1})^{\top},~~\tilde{F}=(\tilde{F}_1,\cdots,\tilde{F}_{N-1})^{\top},\qquad
F_i=(f,\phi_i),~~\tilde{F}_i=(f,\tilde{\phi}_i) .
\nn
\eey
Let 
\bey
&& V=(V_1,\cdots, V_{N-1})^\top,\qquad V_i=u_h(x_i)-\hat{u}_h(x_i)=u_h(x_i)-u_{\tilde{h}}(\tilde{x}_i),
\nn \\
&& v=u_h-\hat{u}_h=\sum\limits_{i=1}^{N-1}V_i\phi_i(x)\in V_h .
\nn
\eey
By subtracting the second equation from the first one in (\ref{lem5.5-2}), re-grouping the terms, 
and taking the inner product of the resulting equation with $V$, we obtain
\beq
V^{\top} AV+V^{\top} (A-\tilde{A})\tilde{U}=V^{\top} (F-\tilde{F}) .
\label{lem5.5-3}
\eeq

We now estimate the terms in (\ref{lem5.5-3}) separately.
First, from Lemma \ref{lem3.2} we have
\beq
 V^{\top} AV = B(v,v)\ge a_0\|v'\|_\Omega^2 .
\label{lem5.5-4}
\eeq
It is not difficult to verify that
\beq
\|v'\|_\Omega^2 = \sum^N_{i=1}\frac{1}{h_i}|V_i-V_{i-1}|^2 .
\label{lem5.5-5}
\eeq
Moreover, from Poincare's inequality we have
\bey
\|v'\|_\Omega^2 & \ge & 8 \| v \|_\Omega^2 = 8 \| \sum_i V_i \phi_i \|_\Omega^2
\nn \\
& = & 8 \sum_{i=1}^{N-1} V_i
\left [\frac{(h_i+h_{i+1})}{3} V_i + \frac{h_i}{6} V_{i-1} + \frac{h_{i+1}}{6} V_{i+1}\right ]
\nn \\
& \ge & \frac{8}{6}  \sum^{N-1}_{i=1}(h_i+h_{i+1})V_i^2 .
\label{lem5.5-6}
\eey
Combining (\ref{lem5.5-4})--(\ref{lem5.5-6}) we get
\bey
&& V^{\top} A V
\nn \\
&\ge& \frac{a_0}{4} \| v' \|_\Omega^2 + \frac{a_0}{8} \| v' \|_\Omega^2
+ \frac{a_0}{4} \| v' \|_\Omega^2 + \frac{3 a_0}{8} \| v' \|_\Omega^2
\nn \\
&\ge & \frac{a_0}4\|u'_h-\hat{u}'_h\|^2_\Omega+a_0\|u_h-\hat{u}_{h}\|^2_\Omega
+ \frac{a_0}{4}\sum^N_{i=1}\frac{1}{h_i}|V_i-V_{i-1}|^2
+\frac{a_0}{2}\sum^N_{i=1} (h_i+h_{i+1}) V_i^2.
\label{lem5.5-7}
\eey

Next, we estimate the term $V^{\top} (F-\tilde{F})$. Noticing that $\phi_{i} + \phi_{i+1} = 1$ on
$(x_{i}, x_{i+1})$ and $\tilde{\phi}_{i} + \tilde{\phi}_{i+1} = 1$ on $(\tilde{x}_{i}, \tilde{x}_{i+1})$, we have
\bey
F_i-\tilde{F}_i  &=& 
\int^{x_{i+1}}_{x_{i-1}}f \phi_idx-\int^{\tilde{x}_{i+1}}_{\tilde{x}_{i-1}}f\tilde{\phi}_idx
\no\\
& = & \left ( \int^{x_{i}}_{x_{i-1}}f
\phi_idx-\int^{\tilde{x}_{i}}_{\tilde{x}_{i-1}}f\tilde{\phi}_idx \right )
+ \left (\int^{x_{i+1}}_{x_{i}}f \phi_idx-\int^{\tilde{x}_{i+1}}_{\tilde{x}_{i}}f\tilde{\phi}_idx\right )
\nn \\
&=& \left (\int^{x_{i}}_{x_{i-1}}f
\phi_idx-\int^{\tilde{x}_{i}}_{\tilde{x}_{i-1}}f\tilde{\phi}_idx \right )
- \left ( \int^{x_{i+1}}_{x_{i}}f \phi_{i+1}dx-\int^{\tilde{x}_{i+1}}_{\tilde{x}_{i}}f\tilde{\phi}_{i+1}dx\right )
\nn \\
&& \quad 
+ \; \left (\int^{x_{i+1}}_{x_{i}}f dx-\int^{\tilde{x}_{i+1}}_{\tilde{x}_{i}}f dx \right )
\nn \\
&=& \left (\int^{x_{i}}_{x_{i-1}}f
\phi_idx-\int^{\tilde{x}_{i}}_{\tilde{x}_{i-1}}f\tilde{\phi}_idx \right )
- \left ( \int^{x_{i+1}}_{x_{i}}f \phi_{i+1}dx-\int^{\tilde{x}_{i+1}}_{\tilde{x}_{i}}f\tilde{\phi}_{i+1}dx\right )
\nn \\
&& \quad 
+ \; \left (\int^{x_{i+1}}_{\tilde{x}_{i+1}}f dx-\int^{x_i}_{\tilde{x}_{i}}f dx \right ) .
\nn
\eey
Thus,
\beq
V^{\top} (F-\tilde{F}) = \sum_{i=1}^N
\left (\int^{x_{i}}_{x_{i-1}}f
\phi_idx-\int^{\tilde{x}_{i}}_{\tilde{x}_{i-1}}f\tilde{\phi}_idx \right ) (V_i-V_{i-1})
+ \sum_{i=1}^{N} (V_{i-1}-V_{i}) \int^{x_i}_{\tilde{x}_{i}}f dx .
\label{lem5.5-8}
\eeq
Denote
\beq
x_i^-=\min\{x_i,\tilde{x}_{i}\},~~x_i^+=\max\{x_i,\tilde{x}_{i}\}.
\label{lem5.5-26}
\eeq
When $(x_{i-1}, x_{i})$ and $(\tilde{x}_{i-1}, \tilde{x}_{i})$ overlap, we have
\bey
&& \left | \int^{x_{i}}_{x_{i-1}}f \phi_idx-\int^{\tilde{x}_{i}}_{\tilde{x}_{i-1}}f\tilde{\phi}_idx\right |
\nn \\
& = &
\left | \int^{x_{i-1}^+}_{x_{i-1}}f \phi_idx + \int^{x_{i}}_{x_{i}^-}f \phi_idx
+ \int^{x_{i}^-}_{x_{i-1}^+}f (\phi_i-\tilde{\phi}_i)dx
-\int^{\tilde{x}_{i}}_{x_{i}^-}f\tilde{\phi}_idx
-\int^{\tilde{x}_{i-1}^+}_{x_{i-1}}f\tilde{\phi}_idx
\right |
\nn \\
& \le & \| f \|_{L^\infty(\Omega)} \left ( |x_{i-1}-\tilde{x}_{i-1}| + |x_{i}-\tilde{x}_{i}| \right )
\nn \\
&& 
+ \| f \|_{L^\infty(\Omega)} \int^{x_{i}^-}_{x_{i-1}^+} \left |\frac{(\tilde{h}_{i}-h_{i})}{h_{i}\tilde{h}_{i}}
(x-x_{i-1})+\frac{(\tilde{x}_{i-1}-x_{i-1})}{\tilde{h}_{i}}\right | d x
\nn \\
& \le & \| f \|_{L^\infty(\Omega)} \left (2 |x_{i-1}-\tilde{x}_{i-1}| +  |x_{i}-\tilde{x}_{i}| + |\tilde{h}_{i}-h_{i} | \right )
\nn \\
& \le & 3 \| f \|_{L^\infty(\Omega)} \left ( |x_{i-1}-\tilde{x}_{i-1}| + |x_{i}-\tilde{x}_{i}| \right ) .
\nn
\eey
On the other hand, if $(x_{i-1}, x_{i})$ and $(\tilde{x}_{i-1}, \tilde{x}_{i})$ do not overlap, we have
\bey
\left | \int^{x_{i}}_{x_{i-1}}f \phi_idx-\int^{\tilde{x}_{i}}_{\tilde{x}_{i-1}}f\tilde{\phi}_idx\right |
& = &
\int^{x_{i-1}^+}_{x_{i-1}} |f \phi_i|dx
+ \int^{\tilde{x}_{i}}_{x_{i}^-}| f\tilde{\phi}_i| dx
\nn \\
& \le & \| f \|_{L^\infty(\Omega)} \left ( |x_{i-1}-\tilde{x}_{i-1}| + |x_{i}-\tilde{x}_{i}| \right ) .
\nn
\eey
For both cases we thus have
\beq
\left | \int^{x_{i}}_{x_{i-1}}f \phi_idx-\int^{\tilde{x}_{i}}_{\tilde{x}_{i-1}}f\tilde{\phi}_idx\right |
\le 3 \| f \|_{L^\infty(\Omega)} \left ( |x_{i-1}-\tilde{x}_{i-1}| + |x_{i}-\tilde{x}_{i}| \right ) .
\label{lem5.5-9}
\eeq
Inserting (\ref{lem5.5-9}) into (\ref{lem5.5-8}) and using Young's inequality and $h_i \le 2/N$,
we get
\beq
|V^{\top} (F-\tilde{F})|  \le \frac{a_0}{20}\sum_{i=1}^{N-1}h_iV_i^2
+C \|X-\tilde{X}\|_\infty^2.
\label{lem5.5-12}
\eeq

We now proceed to estimate $V^{\top} (A-\tilde{A})\tilde{U}$. We start with computing
the non-zero entries of $A = (a_{ij})$. Noticing that
\bey
&& \phi_i \phi_{i-1}' = - \phi'_{i} \phi_{i} = - \frac{1}{2} (\phi_{i}^2)',\quad \phi_i + \phi_{i-1} = 1,\quad
\mbox{ on } (x_{i-1}, x_{i})
\nn \\
&& \phi_{i} \phi_{i+1}' = \phi_{i+1}' - \phi_{i+1}' \phi_{i+1} = \phi_{i+1}'- \frac{1}{2} (\phi_{i+1}^2)',
\quad \phi_i + \phi_{i+1} = 1,\quad
\mbox{ on } (x_{i}, x_{i+1})
\nn
\eey
by direct calculation we have
\bey
a_{i,i-1}&=&B(\phi_{i-1},\phi_{i})
\nn \\
& =& \int^{x_{i}}_{x_{i-1}}\left [a
\phi_{i}'\phi_{i-1}'+b\phi_{i}\phi_{i-1}'+c\phi_{i}\phi_{i-1}\right ]dx\no\\
&=&-\frac{a_i}{h_i}-\frac12b(x_{i})+\frac12\int^{x_{i}}_{x_{i-1}}
b'\phi_{i}dx+\int^{x_{i}}_{x_{i-1}}(c-\frac12b')\phi_{i}\phi_{i-1}dx,
\label{lem5.5-13}
\\
a_{i,i+1}&=&B(\phi_{i+1},\phi_i)
\nn \\
&=&\int^{x_{i+1}}_{x_{i}}\left [a
\phi_{i}'\phi_{i+1}'+b\phi_{i}\phi_{i+1}'+c\phi_{i}\phi_{i+1}\right ]dx\no\\
&=&-\frac{a_{i+1}}{h_{i+1}}+\frac12b(x_{i+1})-\frac12\int^{x_{i+1}}_{x_i}b'\phi_{i+1}dx
+\int^{x_{i+1}}_{x_{i}}(c-\frac12b')\phi_{i}\phi_{i+1}dx,
\label{lem5.5-14}
\\
a_{i,i}&=&B(\phi_i,\phi_{i})
\nn \\
&=&\int^{x_{i+1}}_{x_{i-1}}\left [a
\phi_{i}'\phi_{i}'+(c-\frac12b')\phi_{i}^2\right ]dx\no\\
&=&\frac{a_i}{h_i}+\frac{a_{i+1}}{h_{i+1}}+\int^{x_{i+1}}_{x_{i-1}}(c-\frac12b')\phi_{i}^2dx,
\label{lem5.5-15}
\eey
where
\[
a_i = \frac{1}{h_i} \int_{x_{i-1}}^{x_i} a(x) d x .
\]
Matrix $\tilde{A}=(\tilde{a}_{ij})$ has similar expressions.

Using the above expressions, we have
\bey
& & V^{\top}(A-\tilde{A})\tilde{U}
\nn \\
&=&
\sum^{N-1}_{i=1}V_i\left [(\frac{a_i}{h_i}-\frac{\tilde{a}_i}{\tilde{h}_i}+\frac{a_{i+1}}
 {h_{i+1}}-\frac{\tilde{a}_{i+1}}{\tilde{h}_{i+1}})
 \tilde{U}_{i}-(\frac{a_i}{h_i}-\frac{\tilde{a}_i}{\tilde{h}_i})\tilde{U}_{i-1}-(\frac{a_{i+1}}{h_{i+1}}
 -\frac{\tilde{a}_{i+1}}{\tilde{h}_{i+1}})\tilde{U}_{i+1}\right ]
 \no\\
 &-&\sum^{N-1}_{i=1}\frac12V_i\left [(b(x_{i})-b(\tilde{x}_{i}))\tilde{U}_{i-1}
 -(b(x_{i+1})-b(\tilde{x}_{i+1}))\tilde{U}_{i+1}\right ]
 \no\\
 &+&\sum^{N-1}_{i=1}\frac12V_i\left [(B_{i}-\tilde{B}_{i})\tilde{U}_{i-1}
 +(B_{i+1}-\tilde{B}_{i+1})\tilde{U}_{i+1}\right  ]
 \no\\
&+&\sum^{N-1}_{i=1}V_i\left [(C_{i}-\tilde{C}_{i})\tilde{U}_{i-1}+(C_{2,i}-\tilde{C}_{2,i})
\tilde{U}_i+(C_{i+1}-\tilde{C}_{i+1})\tilde{U}_{i+1}\right ]\no\\
&=&I_1+I_2+I_3+I_4 ,
\label{lem5.5-16}
\eey
where
\[
B_{i} = \int^{x_{i}}_{x_{i-1}} b'\phi_{i}dx, \quad C_{i} = \int^{x_{i}}_{x_{i-1}}(c-\frac12b')\phi_{i}\phi_{i-1}dx,
\quad C_{2,i} = \int^{x_{i+1}}_{x_{i-1}}(c-\frac12b')\phi_{i}^2dx .
\]
Noticing that Lemma \ref{lem5.4} implies
$\| u_{\tilde{h}}' \|_{L^\infty(\Omega)} \le C$, we have
\bey
|I_1| &=&\left|
\sum^{N-1}_{i=1}V_i\left [(\frac{a_i}{h_i}-\frac{\tilde{a}_i}{\tilde{h}_i})(\tilde{U}_i-\tilde{U}_{i-1})
+(\frac{a_{i+1}}{h_{i+1}}-\frac{\tilde{a}_{i+1}}{\tilde{h}_{i+1}})
 (\tilde{U}_{i}-\tilde{U}_{i+1})\right ]\right|
\no\\
&=&\left| \sum^{N}_{i=1}(V_i-V_{i-1})(\frac{a_i}{h_i}
-\frac{\tilde{a}_i}{\tilde{h}_i})(\tilde{U}_i-\tilde{U}_{i-1})\right|\no\\
&\le&C\sum^{N}_{i=1}|V_i-V_{i-1}|~|\frac{a_i}{h_i}-\frac{\tilde{a}_i}{\tilde{h}_i}|\tilde{h}_i
\no\\
&\le&\frac{a_0}{40}\sum^{N}_{i=1}\frac1{h_i}|V_i-V_{i-1}|^2+
C\sum^{N}_{i=1}|a_i\tilde{h}_i-\tilde{a}_ih_i|^2\frac1{h_i} .
\nn 
\eey
From the definitions of $a_i$ and $\tilde{a}_i$,
\bey
a_i\tilde{h}_i-\tilde{a}_ih_i & = & a_i (\tilde{h}_i - h_i) + h_i (a_i - \tilde{a}_i)
\nn \\
& = & a_i (\tilde{h}_i - h_i) + \frac{(\tilde{h}_i - h_i) h_i}{\tilde{h}_i} \cdot
\frac{1}{h_i} \int_{x_{i-1}}^{x_i} a(x) d x 
+ \frac{h_i}{\tilde{h}_i} \left ( \int_{x_{i-1}}^{x_i} a(x) d x - \int_{\tilde{x}_{i-1}}^{\tilde{x}_i} a(x) d x\right ) .
\nn
\eey
From this it is not difficult to obtain
\[
|a_i\tilde{h}_i-\tilde{a}_ih_i|^2\frac1{h_i} \le \frac{C h^2}{ \underline{h}^3} ( |x_i-\tilde{x}_i| + |x_{i-1}-\tilde{x}_{i-1}|) .
\]
Thus, from the fact that $1/(\rho_0 N) \le \underline{h}\le  h_i \le 2/N$ it follows
\beq
|I_1| \le \frac{a_0}{40}\sum^{N}_{i=1}\frac{1}{h_i}|V_i-V_{i-1}|^2+
C N^2 \|X-\tilde{X}\|_\infty^2 .
\label{lem5.5-17}
\eeq
Similarly, we have
\bey
|I_2| & = & \frac{1}{2} \left | \sum_{i=1}^{N-1} (b(x_{i})-b(\tilde{x}_i)) (V_i \tilde{U}_{i-1}- V_{i-1} \tilde{U}_{i}) \right |
\nn \\
&\le&\frac12
\sum^{N-1}_{i=1}~\overline{b'} ~|x_{i}-\tilde{x}_{i}|\left (|V_i| \; |\tilde{U}_i-\tilde{U}_{i-1}|
+ |V_i-V_{i-1}| \; |\tilde{U}_i| \right )
\no\\
&\le&\frac{a_0}{20}\sum^{N-1}_{i=1}h_{i}|V_i|^2+C \|X-\tilde{X}\|_\infty^2
+ \frac{a_0}{40} \sum_{i=1}^N \frac{1}{h_i} |V_i-V_{i-1}|^2 .
\label{lem5.5-18}
\eey

The difference $B_{i}-\tilde{B}_{i}$ involved in  $I_3$
can be estimated in the same manner as for $(\int_{x_{i-1}}^{x_i} f \phi dx
- \int_{\tilde{x}_{i-1}}^{\tilde{x}_i} f \tilde{\phi} dx)$
in $(F_i - \tilde{F}_i)$. We have
\beq
|B_{i}-\tilde{B}_{i}| \le 
3~ \overline{b'}~(|x_{i-1}-\tilde{x}_{i-1}|+|x_i-\tilde{x}_i|) .
\label{lem5.5-19}
\eeq
Using these estimates we obtain
\bey
|I_3| & = & \frac{1}{2} \left | \sum_{i=1}^{N-1} (B_i-\tilde{B}_i) (V_i \tilde{U}_{i-1}- V_{i-1} \tilde{U}_{i}) \right |
\nn \\
&\le&\frac12
\sum^{N-1}_{i=1}~\overline{b'} ~|B_i-\tilde{B}_i |\left (|V_i| \; |\tilde{U}_i-\tilde{U}_{i-1}|
+ |V_i-V_{i-1}| \; |\tilde{U}_i| \right )
\nn \\
&\le&\frac{a_0}{20}\sum^{N-1}_{i=1}h_{i}|V_i|^2+C \|X-\tilde{X}\|_\infty^2
+ \frac{a_0}{40} \sum_{i=1}^N \frac{1}{h_i} |V_i-V_{i-1}|^2 .
\label{lem5.5-20}
\eey
To estimate $I_4$, we denote
\[
C_{0,i} = \int^{x_{i}}_{x_{i-1}}(c-\frac12b') d x,\quad 
C_{1,i} = \int^{x_{i}}_{x_{i-1}}(c-\frac12b') \phi_{i} d x .
\]
Like for $B_{i}-\tilde{B}_{i}$, we have the estimates 
\bey
|C_{i}-\tilde{C}_{i}| & \le & 
\overline{(c-b'/2)}(3|x_{i-1}-\tilde{x}_{i-1}|+3|x_i-\tilde{x}_i|),
\label{lem5.5-21}
\\
|C_{1,i}-\tilde{C}_{1,i}| & \le &
\overline{(c-b'/2)}(3|x_{i-1}-\tilde{x}_{i-1}|+3|x_i-\tilde{x}_i|).
\label{lem5.5-22}
\eey
Moreover,
\bey
C_{2,i} & = & \int^{x_{i}}_{x_{i-1}}(c-\frac12b')\phi_{i}^2dx
+ \int^{x_{i+1}}_{x_{i}}(c-\frac12b')\phi_{i}^2dx
\nn \\
& = &  \int^{x_{i}}_{x_{i-1}}(c-\frac12b')\phi_{i} (1 - \phi_{i-1}) dx
+ \int^{x_{i+1}}_{x_{i}}(c-\frac12b')\phi_{i} (1 - \phi_{i+1}) dx
\nn \\
& = & C_{1,i} - C_{i} + \int^{x_{i+1}}_{x_{i}}(c-\frac12b')\phi_{i} d x - C_{i+1} 
\nn \\
& = & C_{1,i} - C_{i} + \int^{x_{i+1}}_{x_{i}}(c-\frac12b')(1-\phi_{i+1}) d x - C_{i+1} 
\nn \\
& = & C_{0,i+1} + C_{1,i}- C_{1,i+1} - C_{i}  - C_{i+1} .
\nn
\eey
Then,
\bey
I_4 & = & \sum_{i=1}^{N-1} V_i \left [ (C_{i}-\tilde{C}_{i}) (\tilde{U}_{i-1}-\tilde{U}_{i})
- (C_{i+1}-\tilde{C}_{i+1}) (\tilde{U}_{i}-\tilde{U}_{i+1}) \right .
\nn \\
&& \quad \left .
+ ((C_{1,i}-\tilde{C}_{1,i})-(C_{1,i+1}-\tilde{C}_{1,i+1})) \tilde{U}_i
 + (C_{0,i+1}-\tilde{C}_{0,i+1}) \tilde{U}_i \right ]
 \nn \\
 & = & \sum_{i=1}^{N-1} (V_i-V_{i-1}) (C_{i}-\tilde{C}_{i}) (\tilde{U}_{i-1}-\tilde{U}_{i})
 + \sum_{i=1}^{N-1} (C_{1,i}-\tilde{C}_{1,i}) (V_i \tilde{U}_{i}-V_{i-1} \tilde{U}_{i-1})
 \nn \\
 && \quad + \sum_{i=1}^{N-1} (V_i \tilde{U}_{i}-V_{i-1} \tilde{U}_{i-1}) \int_{\tilde{x}_i}^{x_i} (c-\frac{1}{2} b')  d x .
 \nn
\eey
From this we get
\bey
|I_4|&\le& \frac{a_0}{20}\sum_{i=1}^{N-1}(h_i+h_{i+1})V_i^2+C \|X-\tilde{X}\|_\infty^2
+ \frac{a_0}{40} \sum_{i=1}^N \frac{1}{h_i} |V_i-V_{i-1}|^2 .
\label{lem5.5-24}
\eey
Combining (\ref{lem5.5-16}), (\ref{lem5.5-17}), (\ref{lem5.5-18}), (\ref{lem5.5-20}),
and (\ref{lem5.5-24}), we get
\beq
V^{\top}(A-\tilde{A})\tilde{U} \le
\frac{a_0}8\sum^{N}_{i=1}\frac{1}{h_i}|V_i-V_{i-1}|^2 
+ \frac{a_0}{4}\sum_{i=1}^{N-1}(h_i+h_{i+1})V_i^2
+ C N^2 \|X-\tilde{X}\|_\infty^2 .
\label{lem5.5-25}
\eeq

Finally, (\ref{lem5.5-1}) follows from (\ref{lem5.5-7}), (\ref{lem5.5-12}), and (\ref{lem5.5-25}).
\proofend

\vspace{10pt}

{\bf Proof of Theorem \ref{thm2.3}.}  From Lemma \ref{lem5.5}
we can see that the key to the proof of this theorem
is to estimate $\|\hat{u}'_h-u'_{\tilde{h}}\|_{L^1(\Omega)}$ and
$\|\hat{u}_h-u_{\tilde{h}}\|_{L^1(\Omega)}$.
For this purpose, we notice from assumption (\ref{thm2.3-1}) that
$\| X - \tilde{X}\|_\infty < \min_i\{h_i, \tilde{h}_i\}$ and 
\[
(x_{i-1}, x_{i}) \cap (\tilde{x}_{i-1}, \tilde{x}_{i}) \neq \emptyset,\qquad i = 1, ..., N.
\]
As a consequence, we can divide $[x_{i-1},x_i]$
into subintervals $[x_{i-1},x^+_{i-1}]$, $[x_{i-1}^+,x_i^-]$, and $[x^-_i,x_i]$. On these intervals
$\hat{u}_h-u_{\tilde{h}}$ can be expressed as
\bey
(\hat{u}_h-u_{\tilde{h}})|_{[x_{i-1},x^+_{i-1}]} &=&
\left (u_{\tilde{h}}(\tilde{x}_i)\phi_i(x)+u_{\tilde{h}}(\tilde{x}_{i-1})\phi_{i-1}(x)\right )
-\left (u_{\tilde{h}}(\tilde{x}_{i-1})\tilde{\phi}_{i-1}(x)+u_{\tilde{h}}(\tilde{x}_{i-2})\tilde{\phi}_{i-2}(x)\right )
\no\\
&=&\left (u_{\tilde{h}}(\tilde{x}_{i})-u_{\tilde{h}}(\tilde{x}_{i-1})\right )\phi_{i}(x)
+\left (u_{\tilde{h}}(\tilde{x}_{i-1})-u_{\tilde{h}}(\tilde{x}_{i-2})\right )\phi_{i-1}(x)
\no\\
&&\qquad + \left (u_{\tilde{h}}(\tilde{x}_{i-1})-u_{\tilde{h}}(\tilde{x}_{i-2})\right )
\left (\phi_{i}(x)-\tilde{\phi}_{i-1}(x)\right ),
\label{thm2.3-5}
\\
(\hat{u}_h-u_{\tilde{h}})|_{[x^+_{i-1},x^-_{i}]} &=&
\left (u_{\tilde{h}}(\tilde{x}_{i})-u_{\tilde{h}}(\tilde{x}_{i-1})\right )\left (\phi_{i}(x)-\tilde{\phi}_{i}(x)\right ),
\label{thm2.3-6}
\\
(\hat{u}_h-u_{\tilde{h}})|_{[x_i^-,x_{i}]} &=&
\left (u_{\tilde{h}}(\tilde{x}_{i})-u_{\tilde{h}}(\tilde{x}_{i+1})\right )\phi_{i}(x)
+\left (u_{\tilde{h}}(\tilde{x}_{i-1})-u_{\tilde{h}}(\tilde{x}_{i})\right )\phi_{i-1}(x)
\no\\
&&\qquad + \left (u_{\tilde{h}}(\tilde{x}_{i+1})-u_{\tilde{h}}(\tilde{x}_{i})\right )
\left (\phi_{i}(x)-\tilde{\phi}_{i+1}(x)\right ).
\label{thm2.3-7}
\eey
Integrating $|\hat{u}_h-u_{\tilde{h}}|$ over the subintervals and 
using the above expressions and Lemma \ref{lem5.4}, we get
\bey
\|\hat{u}_h-u_{\tilde{h}}\|_{L^1(x_{i-1},x^+_{i-1})} &\le&
C(\tilde{h}_i+\tilde{h}_{i-1})|x_{i-1}-\tilde{x}_{i-1}|,
\no\\
\|\hat{u}_h-u_{\tilde{h}}\|_{L^1(x^+_{i-1},x^-_{i})} &\le&
Ch_i(|x_i-\tilde{x}_{i}|+|x_{i-1}-\tilde{x}_{i-1}|),
\no\\
\|\hat{u}_h-u_{\tilde{h}}\|_{L^1(x^-_{i},x_{i})} &\le&
C(\tilde{h}_i+\tilde{h}_{i+1})|x_i-\tilde{x}_{i}|.
\no
\eey
Summing these estimates from $i=1$ to $i=N$ yields
\bey
\|\hat{u}_h-u_{\tilde{h}}\|_{L^1(\Omega)} &\le&
C\|X-\tilde{X}\|_\infty.
\label{thm2.3-8}
\eey

Moreover, differentiating (\ref{thm2.3-5})-(\ref{thm2.3-7}) leads to
\bey
(\hat{u}'_h-u'_{\tilde{h}})|_{[x_{i-1},x^+_{i-1}]} &=&
u'_{\tilde{h}}(\tilde{x}_{i})\frac{\tilde{h}_i}{h_i}-u'_{\tilde{h}}(\tilde{x}_{i-1}),
\label{thm2.3-9}
\\
(\hat{u}'_h-u'_{\tilde{h}})|_{[x^+_{i-1},x^-_{i}]} &=&
u'_{\tilde{h}}(\tilde{x}_{i})(\frac{\tilde{h}_i}{h_i}-1),
\label{thm2.3-10}
\\
(\hat{u}'_h-u'_{\tilde{h}})|_{[x_i^-,x_{i}]}
&=&u'_{\tilde{h}}(\tilde{x}_{i})\frac{\tilde{h}_i}{h_i}
-u'_{\tilde{h}}(\tilde{x}_{i+1}).
\label{thm2.3-11}
\eey
Integrating $|\hat{u}'_h-u'_{\tilde{h}}|$ over the subintervals and 
using Lemma \ref{lem5.4}, we obtain
\bey
\|\hat{u}'_h-u'_{\tilde{h}}\|_{L^1(x_{i-1},x^+_{i-1})} &\le&
C(1+\frac{\tilde{h}_i}{h_i})|x_{i-1}-\tilde{x}_{i-1}|,
\label{thm2.3-12}
\\
\|\hat{u}'_h-u'_{\tilde{h}}\|_{L^1(x^+_{i-1},x^-_{i})} &\le&
C|h_i-\tilde{h}_i|,
\label{thm2.3-13}
\\
\|\hat{u}'_h-u'_{\tilde{h}}\|_{L^1(x_i^-,x_{i})}
&\le&C(1+\frac{\tilde{h}_i}{h_i})|x_i-\tilde{x}_{i}|.
\label{thm2.3-14}
\eey
Thus, combining these estimates gives
\beq
\|(\hat{u}_h-u_{\tilde{h}})'\|_{L^1(\Omega)}\le C
N\|X-\tilde{X}\|_\infty.
\label{thm2.3-15}
\eeq
Inequality (\ref{thm2.3-2}) follows from the above estimate, Lemma \ref{lem5.5}, and the triangle inequality.

Next, recalling from Lemma \ref{lem5.5} and (\ref{thm2.3-8}) that
\[
\|u_h-\hat{u}_{h}\|_\Omega\le C~N \|X-\tilde{X}\|_\infty,~~\|\hat{u}_h-u_{\tilde{h}}\|_{L^1(\Omega)}
\le C\|X-\tilde{X}\| ,
\]
by Schwarz' inequality and Lemma \ref{lem5.4} we have
\bey
\|u_h-u_{\tilde{h}}\|_\Omega^2 & = & \int_\Omega | u_h-u_{\tilde{h}} | \cdot | u_h-u_{\tilde{h}} | d x
\nn \\
& \le & \int_\Omega | u_h-u_{\tilde{h}} | \cdot | u_h-\hat{u}_h | d x
+ \int_\Omega | u_h-u_{\tilde{h}} | \cdot | \hat{u}_h-u_{\tilde{h}} | d x
\nn \\
&\le&
\|u_h-u_{\tilde{h}}\|_\Omega\|u_h-\hat{u}_h\|_\Omega
+\|u_h-u_{\tilde{h}}\|_{L^\infty(\Omega)}\|\hat{u}_h-u_{\tilde{h}}\|_{L^1(\Omega)}
\no\\
& \le & C N^{-2} \cdot N \|X-\tilde{X}\|_{\infty} + C N^{-\frac{3}{2}} \cdot \|X-\tilde{X}\|_{\infty}
\no\\
&\le&CN^{-1}\|X-\tilde{X}\|_\infty,\no\eey
 which gives
(\ref{thm2.3-3}).

Finally, the finite element equation (\ref{fem-1}) implies that
\[
\|u_h\|_E^2 - \|u_{\tilde{h}}\|_E^2 = (f,u_h-u_{\tilde{h}}) .
\]
From (\ref{thm2.3-2}) we have
\[
\left | \|u_h\|_E^2 - \|u_{\tilde{h}}\|_E^2 \right | = \|f\|_\Omega \| u_h-u_{\tilde{h}}\|_\Omega
\le C N^{-\frac12}\|X-\tilde{X}\|^{\frac12}_\infty ,
\]
which gives (\ref{thm2.3-4}).
\proofend

% section 6
\section{ Existence of equidistributing meshes}
\label{SEC:existence}

We prove Theorem \ref{thm2.4} in this section.  The existence of equidistributing meshes
is equivalent to the existence of fixed points of the map $G_N$ defined
by the iterative algorithm in \S\ref{SEC:numerical-results}. The key is
to show that $G_N$ maps $S_N$ into $S_N$ and is continuous.

Recall from (\ref{eq-4}) that the mesh $Y=G_NX$ satisfies
\beq
\int^{y_i}_{y_{i-1}}\rho (x)dx=\frac{\s_h}{N}, \quad i = 1, ..., N
\label{eq-5}
\eeq
where $\rho(x)$ and $\sigma_h$ are defined in (\ref{rho-1}) and (\ref{sigma-1}) based on
the solution $u_h$ obtained on mesh $X$ (i.e. $\pi_h$).

\begin{lem}
\label{lem6.0}
Assume that $X \in S_N$ and $u \in H^2(\Omega)$. Then there exists a positive integer $N_0$,
independent of the finite element approximation and the mesh, such that, for any $N \ge N_0$, 
\beq
 \frac12\|r\|^{\frac23}_{L^{\frac23}(\Omega)} \le \a_h^{\frac13}\le
\frac32\|r\|^{\frac23}_{L^{\frac23}(\Omega)}.
\label{lem6.0-1}
\eeq
Moreover, for any $N\ge N_0$,
\beq
1\le \rho \le \rho_0, \quad 1\le \s_h \le 2 ,
\label{lem6.0-2}
\eeq
where $\rho_0$ is a constant defined in (\ref{rho0}).
\end{lem}

{\bf Proof.} The existence of $N_0$ is guaranteed by Lemma \ref{lem3.10}.
Its independence of the finite element approximation and the mesh is clear from (\ref{lem3.10-1})
for the situation $r'\in L^1(\Omega)$. For the situation $r \in L^2(\Omega)$,
we can choose $\epsilon = \theta \| r \|_{L^{\frac 2 3}(\Omega)}$ for some value of $\theta$
(cf. the proof of Theorem \ref{thm2.1}). Then, a smoother function $\tilde{r}$ which
is independent of the approximation and the mesh can be chosen and an inequality similar to 
(\ref{thm2.1-9}) can be obtained. Thus, an $N_0$ independent of the finite element
approximation and the mesh also exists for the situation  $r \in L^2(\Omega)$.

The inequality $\rho_i \ge 1$ follows immediately from the definition of $\rho$. For the upper bound of
$\rho$, from (\ref{rho-2}), (\ref{lem6.0-1}), Lemma \ref{lem5.2}, Young's inequality,
and the inequalities $\|r\|_{L^2(\Omega)} \le \|r\|_{L^\infty(\Omega)}$ and $N \ge 1$ we have
\bey
\rho_{i}  & = &(1+\a_h^{-1}\left <r_h\right >_i^2)^{\frac13}
 \nn \\
 & \le &1+\a_h^{-\frac13}\|r_h\|^{\frac23}_{L^\infty(\Omega)}
\no\\
& = &1+\a_h^{-\frac13}\|r - a' e_h' + b e_h' + c e_h \|_{L^\infty(\Omega)}^{\frac23}
\nn \\
& \le & 
1+ 2\|r\|_{L^{\frac23}(\Omega)}^{-\frac23} \left ( \|r\|_{L^\infty(\Omega)}
+ \overline{(b-a')} \| e_h'\|_{L^\infty(\Omega)}
+\bar{c}\| e_h\|_{L^\infty(\Omega)} \right )^{\frac23}
\nn \\
& \le & 
1+ 2\|r\|_{L^{\frac23}(\Omega)}^{-\frac23} \|r\|_{L^\infty(\Omega)}^{\frac23}
\left ( 1 + \frac{\sqrt{\rho_0+1} C_1}{\sqrt{N}} 
+\frac{C_2}{ N} \right )^{\frac23}
\nn \\
& \le & 
1+ 2\|r\|_{L^{\frac23}(\Omega)}^{-\frac23} \|r\|_{L^\infty(\Omega)}^{\frac23}
\left ( 1 + C_1+ C_2 + C_1 \sqrt{\rho_0} \right )^{\frac23}
\nn \\
& \le & 
1+ 2\|r\|_{L^{\frac23}(\Omega)}^{-\frac23} \|r\|_{L^\infty(\Omega)}^{\frac23}
\left ( 1 + C_1+ C_2  \right )^{\frac23}
+ 2\|r\|_{L^{\frac23}(\Omega)}^{-\frac23} \|r\|_{L^\infty(\Omega)}^{\frac23}
C_1^{\frac 2 3} \rho_0^{\frac 1 3}
\nn \\
& \le & 
1+ 2\|r\|_{L^{\frac23}(\Omega)}^{-\frac23} \|r\|_{L^\infty(\Omega)}^{\frac23}
\left ( 1 + C_1+ C_2  \right )^{\frac23}
+ 2\|r\|_{L^{\frac23}(\Omega)}^{-\frac23} \|r\|_{L^\infty(\Omega)}^{\frac23}
\left ( 1 + C_1+ C_2  \right )^{\frac 2 3} \rho_0^{\frac 1 2} ,
\nn
\eey
where $C_1$ and $C_2$ denote the constants in (\ref{lem5.2-1}) and (\ref{lem5.2-2}), respectively.
Notice that $\rho_0 \ge 1$ and thus $\rho_0^{\frac 13} \le \rho_0^{\frac 1 2}$.
Letting $\gamma = 2 \left ( 1 + C_1+ C_2  \right )^{\frac 2 3}$, from the definition of $\rho_0$, (\ref{rho0}),
we thus have
\bey
\rho_i
& \le & 
1+ \gamma \|r\|_{L^{\frac23}(\Omega)}^{-\frac23} \|r\|_{L^\infty(\Omega)}^{\frac23}
+ \gamma \|r\|_{L^{\frac23}(\Omega)}^{-\frac23} \|r\|_{L^\infty(\Omega)}^{\frac23} \rho_0^{\frac 1 2}
\nn \\
& = &
\left [1+ \gamma \|r\|_{L^{\frac23}(\Omega)}^{-\frac23} \|r\|_{L^\infty(\Omega)}^{\frac23} \right ]^2
\nn \\
& = & \rho_0 .
\nn
\eey
The bounds for $\sigma_h$ follow from the bounds for $\rho$ and the definitions of $\sigma_h$ and $\alpha_h$.
\proofend

\begin{lem}
\label{lem6.1}
Assume that $u \in H^2(\Omega)$. For any $N \ge N_0$ where $N_0$ is defined
in Lemma \ref{lem6.0}, then $G_N(S_N) \subset S_N$ or $G_N: S_N \to S_N$. 
\end{lem}

{\bf Proof.} Lemma \ref{lem6.0} implies that for any given $X\in S_N$,
$1\le\rho \le \rho_0$ and $1 \le \sigma_h \le 2$. From (\ref{eq-5}) we then have
\[
(y_i-y_{i-1}) \le  \int^{y_i}_{y_{i-1}}\rho (x)dx=\frac{\s_h}{N}\le \frac2N
\]
and 
\[
(y_i-y_{i-1}) \rho_0 \ge \int^{y_i}_{y_{i-1}}\rho (x)dx = \frac{\s_h}{N}\ge \frac1N.
\]
Thus, $Y = G_N X \in S_N$ for any $X \in S_N$.
\proofend

\begin{lem}
\label{lem6.2}
Assume that $f \in L^\infty(\Omega)$, $u \in H^2(\Omega)$, and $X,\,\tilde{X} \in S_N$.
Then, 
\beq
 \tilde{h}_i\left |\left <r_h\right >^2_{i}-\left <r_{\tilde{h}}\right >^2_{i}\right | \le
C \left (\|X-\tilde{X}\|_{\infty}+ \|u'_h-u'_{\tilde{h}}\|_{L^1(K_i)}+\|u_h-u_{\tilde{h}}\|_{L^1(K_i)}\right ).
\label{lem6.2-1}
\eeq

\end{lem}
{\bf Proof.} Using the notation (\ref{lem5.5-26}) we have
 \bey
 \left <r_h\right >^2_{i}-\left <r_{\tilde{h}}\right >^2_{i}
 &=&\frac{\tilde{h}_i-h_i}{\tilde{h}_i}\left <r_h\right >^2_i
 +\frac1{\tilde{h}_i}\int^{x_i^-}_{x_{i-1}^+}(|r_h|^2-|r_{\tilde{h}}|^2)dx
 \no\\
 &+&\frac1{\tilde{h}_i}(\int^{x_i}_{x_{i}^-}+\int^{x_{i-1}^+}_{x_{i-1}})|r_h|^2dx
 -\frac1{\tilde{h}_i}(\int^{\tilde{x}_i}_{x_{i}^-}+\int^{x_{i-1}^+}_{\tilde{x}_{i-1}})|r_{\tilde{h}}|^2dx.
 \nn
 \eey
 From the assumption $f \in L^\infty(\Omega)$ and Lemma \ref{lem5.4}, we have
 $\| r_h \|_{L^\infty(\Omega)} \le C$ and $\| r_{\tilde{h}} \|_{L^\infty(\Omega)} \le C$.
 It follows that
 \bey
&&\left | \left <r_h\right >^2_{i}-\left <r_{\tilde{h}}\right >^2_{i}\right |
\nn \\
& \le &\frac{1}{\tilde{h}_i}|h_i-\tilde{h}_i|\cdot \|r_h\|^2_{L^\infty(\Omega)}
 + \frac1{\tilde{h}_i}\int^{x_i^-}_{x_{i-1}^+}|r_h-r_{\tilde{h}}|\cdot |r_h+r_{\tilde{h}}|dx
 \no\\
 & & \quad +\frac1{\tilde{h}_i}
 (|x_i-\tilde{x}_i|+|x_{i-1}-\tilde{x}_{i-1}|)(\|r_h\|^2_{L^\infty(\Omega)}+\|r_{\tilde{h}}\|^2_{L^\infty(\Omega)})\no\\
  &\le&\frac{C}{ \tilde{h}_i}\|X-\tilde{X}\|_\infty+\frac{C}{\tilde{h}_i}\|r_h-r_{\tilde{h}}\|_{L^1(K_i)}\|r_h+r_{\tilde{h}}\|_{L^\infty(K_i)}\no\\
  &\le&
  \frac{C}{\tilde{h}_i}\|X-\tilde{X}\|_\infty+\frac{C}{\tilde{h}_i}
  (\|u'_h-u'_{\tilde{h}}\|_{L^1(K_i)}+\|u_h-u_{\tilde{h}}\|_{L^1(K_i)}),\no\eey
which gives (\ref{lem6.2-1}).
\proofend

\begin{lem}
\label{lem6.3}
Assume that $f \in L^\infty(\Omega)$, $u \in H^2(\Omega)$, and $X,\,\tilde{X} \in S_N$ satisfying
$\| X - \tilde{X}\|_\infty < 1/(\rho_0 N)$. We also assume that $N \ge N_0$ where $N_0$ is defined
in Lemma \ref{lem6.0}. Then,
\beq
\int^1_0|\rho -\tilde{\rho}|dx\le  C N\|X-\tilde{X}\|_{\infty}+C(N\|X-\tilde{X}\|_{\infty})^{\frac13}.
\label{lem6.3-1}
\eeq
\end{lem}

{\bf Proof.} Using Lemmas \ref{lem6.0}, \ref{thm2.3}, and \ref{lem6.2} and the inequality
$\|u_h-u_{\tilde{h}}\|_{L^1(\Omega)}\le\|(u_h-u_{\tilde{h}})'\|_{L^1(\Omega)}$,
we have
 \bey
|\a^{\frac13}_h-\a^{\frac13}_{\tilde{h}}| &\le&
\sum^{N}_{i=1}\left [|h_i-\tilde{h}_i|\left <r_h\right >^{\frac23}_i
+\tilde{h}_{i}|\left <r_h\right >^{\frac23}_i-\left <r_{\tilde{h}}\right >^{\frac23}_i|\right ]
\no\\
&\le&C~N\|X-\tilde{X}\|_\infty+\sum^{N}_{i=1}\tilde{h}_i|\left <r_h\right >^2_i-\left <r_{\tilde{h}}\right >^2_i|^{\frac13}
\no\\
&\le&C~N\|X-\tilde{X}\|_\infty+(\sum^{N}_{i=1}\tilde{h}_i|\left <r_h\right >^2_i-\left <r_{\tilde{h}}\right >^2_i|)^{\frac13}
\no\\
&\le&C~N\|X-\tilde{X}\|_\infty+C(\|(u_h-u_{\tilde{h}})'\|_{L^1(\Omega)}+\|u_h-u_{\tilde{h}}\|_{L^1(\Omega)}))^{\frac13}
\no\\
&\le&C~N\|X-\tilde{X}\|_{\infty}+C(N\|X-X'\|_{\infty})^{\frac13},
\label{lem6.3-2}
\\
|\rho_{i}-\tilde{\rho}_{i}|&=&\frac{|\rho^3_{i}-\tilde{\rho}_{i}^3|}{\rho_{i}^2+\rho_{i}\tilde{\rho}_{i}+\tilde{\rho}_{i}^2}
 \no\\
 &\le&
 \frac13 |\frac{\left <r_h\right >_i^{2}}{\a_h}-\frac{\left <r_{\tilde{h}}\right >_i^2}{\a_{\tilde{h}}}|
 \no\\
&\le&
\frac13|\frac{\a_h-\a_{\tilde{h}}}{\a_h\a_{\tilde{h}}}|\left <r_h\right >_i^2+\frac1{3\a_{\tilde{h}}}
|\left <r_h\right >^2_i-\left <r_{\tilde{h}}\right >_i^2|
\no\\
&\le&
\frac13|\frac{|\a^{\frac13}_h-\a^{\frac13}_{\tilde{h}}|}{\a_h\a_{\tilde{h}}}(\a_h^{\frac23}
+\a_h^{\frac13}\a_{\tilde{h}}^{\frac13}+\a_{\tilde{h}}^{\frac23})|\left <r_h\right >_i^2+\frac1{3\a_{\tilde{h}}}
|\left <r_h\right >^2_i-\left <r_{\tilde{h}}\right >_i^2|
\no\\
&\le&\frac{|\a^{\frac13}_h-\a^{\frac13}_{\tilde{h}}|}{\left (\frac 1 2 \|r\|_{L^{\frac 2 3}(\Omega)}^{\frac 2 3}\right )^4}
\left <r_h\right >_i^2
+\frac1{3\a_{\tilde{h}}}|\left <r_h\right >^2_i-\left <r_{\tilde{h}}\right >_i^2|
\no\\
&\le&C|\a^{\frac13}_h-\a^{\frac13}_{\tilde{h}}|
+\frac{C}{\tilde{h}_i}\left (\|X-\tilde{X}\|_\infty+\|u'_h-u'_{\tilde{h}}\|_{L^1(K_i)}+\|u_h-u_{\tilde{h}}\|_{L^1(K_i)}\right ).
\label{lem6.3-3}
\eey
It follows from (\ref{lem6.3-3}), Lemma \ref{lem6.2}, and Theorem {\ref{thm2.3} that
\bey
\int^1_0|\rho -\tilde{\rho}|dx
 &=&
\sum^N_{i=1}(\int^{x_i}_{x_{i}^-}+\int^{x^+_{i-1}}_{x_{i-1}})|\rho -\tilde{\rho}| d x
+\sum^{N}_{i=1}(x_i^{-}-x_{i-1}^+)|\rho_{i}-\tilde{\rho}_{i}|\no\\
 &\le&
\sum^N_{i=1}[(x_i-x_{i}^-)+(x^+_{i-1}-x_{i-1})]\rho_0+\sum^N_{i=1}\tilde{h}_i|\rho_{i}-\tilde{\rho}_{i}|
\no\\
 &\le&C|\a^{\frac13}_h-\a^{\frac13}_{\tilde{h}}|+C N\|X-\tilde{X}\|_{\infty}
 +C(\|(u_h-u_{\tilde{h}})'\|_{L^1(\Omega)}+\|u_h-u_{\tilde{h}}\|_{L^1(\Omega)})
 \no\\
 &\le&C|\a^{\frac13}_h-\a^{\frac13}_{\tilde{h}}|+C N\|X-\tilde{X}\|_{\infty},
 \label{lem6.3-4}
 \eey
 which, together with (\ref{lem6.3-2}), gives (\ref{lem6.3-1}).
 \proofend

\begin{lem}
\label{lem6.4}
Under the assumptions of Lemma \ref{lem6.3}, $G_N$ is
a continuous map from $S_N$ into $S_N$:
\beq
\|G_N X-G_N\tilde{X}\|_\infty\leq
C N\|X-\tilde{X}\|_{\infty}+C(N\|X-\tilde{X}\|_{\infty})^{\frac13}.
\label{lem6.4-1}
\eeq
\end{lem}

{\bf Proof.} Let $Y=G_N X$ and $\tilde{Y}=G_N\tilde{X}$.
From the equidistribution relation (\ref{eq-5}), we obtain
\[
\int_0^{y_i} \rho d x - \int_0^{\tilde{y}_i} \tilde{\rho} d x = \frac{i}{N} \left (\sigma_h - \sigma_{\tilde{h}}\right )
\]
or
\[
\int_{\tilde{y}_i}^{y_i} \rho d x = \frac{i}{N} \int_0^1 (\rho - \tilde{\rho}) d x
+ \int_0^{\tilde{y}_i} (\tilde{\rho}-\rho) d x .
\]
It follows from Lemmas \ref{lem6.0} and \ref{lem6.3} that
\bey
| y_i - \tilde{y}_i| & \le & | \int_{\tilde{y}_i}^{y_i} \rho d x |
\nn \\
& \le & 2 \int_0^1 |\rho - \tilde{\rho}| d x 
\label{lem6.4-2}
\\
&\le & C N\|X-\tilde{X}\|_{\infty}+C(N\|X-\tilde{X}\|_{\infty})^{\frac13},
\nn
\eey
which gives (\ref{lem6.4-1}).
\proofend

\vspace{10pt}

The term involving $(N\|X-\tilde{X}\|_{\infty})^{\frac13}$ in the above lemma can be dropped
for a situation shown in the following lemma.

\begin{lem}
\label{lem6.5}
Assume that the assumptions of Lemma \ref{lem6.3} hold. If further
$r_0=\min_{x\in \Omega}|r(x)|>0$ and $N$ is sufficiently large,
then $G_N: S_N \to S_N$ is a continuous map satisfying
\beq
\|G_N X-G_N \tilde{X}\|_\infty\leq C~N\|X-\tilde{X}\|_\infty.
\label{lem6.5-1}
\eeq
\end{lem}

{\bf Proof.}
Recall from Lemma \ref{lem5.4} that
\[
 \|e_h\|_{L^\infty(\Omega)}\le \frac{C}{N^{\frac32}},\qquad
 \|e'_h\|_{L^\infty(\Omega)}\le\frac{C}{\sqrt{N}}.
\]
There exists a sufficiently large $N$ such that
\bey
|r_h(x)|&=&|f(x)-(b-a')u'_h(x)-cu_h(x)|=|r(x)+(b-a')e_h(x)+ce_h(x)|\no\\
&\ge&
r_0-\overline{(b-a')}\|e'_h\|_{L^\infty(\Omega)}-\bar{c}\|e_h\|_{L^\infty(\Omega)}\no\\
&\ge&\frac12r_0,\qquad \forall x\in \Omega.
\label{lem6.5-2}
\eey
From this and Lemma \ref{lem6.2}, we can estimate $|\a_h^{\frac13}-\a_{\tilde{h}}^{\frac13}|$ as
\bey
|\a^{\frac13}_h-\a^{\frac13}_{\tilde{h}}| &\le&
\sum^{N}_{i=1}\left [|h_i-\tilde{h}_i|\left <r_h\right >^{\frac23}_i
+\tilde{h}_{i}|\left <r_h\right >^{\frac23}_i-\left <r_{\tilde{h}}\right >^{\frac23}_i|\right ]
\no\\
&\le&C~N\|X-\tilde{X}\|_\infty+\sum^{N}_{i=1}\tilde{h}_i
\frac{|\left <r_h\right >^2_i-\left <r_{\tilde{h}}\right >^2_i|}{r^{\frac43}_h(\xi_i)
+r^{\frac23}_h(\xi_i)r^{\frac23}_{\tilde{h}}(\tilde{\xi}_i)+r^{\frac43}_{\tilde{h}}(\tilde{\xi}_i)}
\no\\
&\le&C~N\|X-\tilde{X}\|_\infty
+\sum^{N}_{i=1}\tilde{h}_i\frac{|\left <r_h\right >^2_i-\left <r_{\tilde{h}}\right >^2_i|}{r_0^{\frac43}}
\no\\
&\le& C~N\|X-\tilde{X}\|_\infty,
\label{lem6.5-3}
\eey
where $\xi_i\in K_i,~~\tilde{\xi}_i\in \tilde{K}_i$. Combining (\ref{lem6.5-3}) with
(\ref{lem6.3-4}) gives
\beq
\int^1_0|\rho -\tilde{\rho}|dx
 \le C|\a_h^{\frac13}-\a_{\tilde{h}}^{\frac13}|+C~N\|X-\tilde{X}\|_{\infty}\le C~N\|X-\tilde{X}\|_{\infty}.
 \label{lem6.5-4}
\eeq

 Finally, by (\ref{lem6.5-4})  and (\ref{lem6.4-2}) we obtain
\[
\|Y-\tilde{Y}\|_\infty\leq
2 \int^1_0|\rho -\tilde{\rho}|dx\le
CN\|X-\tilde{X}\|_{\infty},
\]
which gives (\ref{lem6.5-1}).
\proofend

\vspace{10pt}

{\bf Proof of Theorem \ref{thm2.4}.} 
From Lemmas \ref{lem6.1} and \ref{lem6.4} we see that $G_N$ is a continuous map
from $S_N$ to $S_N$. Recall that $S_N$ is a closed, convex set.
By Brouwer's theorem, $G_N$ has at least a fixed point in $S_N$. Since any fixed point of $G_N$ is
an equidistributing mesh, we have proven that an equidistributing mesh exists and is in $S_N$.
\proofend

\vspace{20pt}

{\bf Acknowledgment.} The work was supported in part by the NSF (USA)
under grants DMS-0410545 and DMS-0712935, by the NSF of China under grant 10671154,
and by the National Basic Research Program (China) under grant 2005CB321703.
The work was done while Y. He was visiting
the Department of Mathematics of the University of Kansas from January to July, 2008.

%---------------------------------------------------
%\bibliographystyle{plain}
%\bibliography{/Users/huang/tex/bib/mmesh}

\end{document}